\documentclass[english,11pt]{article}

\usepackage[german,french,english]{babel}
\usepackage[cp850]{inputenc}
\usepackage{latexsym,graphicx, fancybox}
\usepackage{amssymb, amsmath, amsfonts}
\usepackage{color}
\usepackage{latexsym}

\font\tenmath=msbm10 scaled 1200
\font\sevenmath=msbm7 scaled 1200
\font\fivemath=msbm5 scaled 1200

\newcommand{\vertiii}[1]{{\left\vert\kern-0.25ex\left\vert\kern-0.25ex\left\vert #1 
    \right\vert\kern-0.25ex\right\vert\kern-0.25ex\right\vert}}

\newfam\mathfam \textfont\mathfam=\tenmath
\scriptfont\mathfam=\sevenmath \scriptscriptfont\mathfam=\fivemath
\def\math{\fam\mathfam}
\def\R{{\math R}}
\def\N{{\math N}}
\def\E{{\math E}}

\def\P{{\math P}}

\newtheorem{Theorem}{Theorem}[section]

\newtheorem{Proposition}[Theorem]{Proposition}

\newtheorem{Lemma}[Theorem]{Lemma}
\newtheorem{Corollary}[Theorem]{Corollary}

\newcommand{\cM}{M}
\newcommand{\cW}{W}
\newcommand{\cAW}{AW}

\newfam\mathfam \textfont\mathfam=\tenmath
\scriptfont\mathfam=\sevenmath \scriptscriptfont\mathfam=\fivemath
\def\math{\fam\mathfam}

\def \^#1{\if#1i{\accent"5E\i}\else{\accent"5E#1}\fi}

\def \cqfd{\quad_\Box}

\def \lecx{\le_{cvx}}

\begin{document}
\selectlanguage{english}
\title{\bf Quantization and martingale couplings}% {Quantization-based monotone approximations with respect to the convex order}
 
\author{ 
{\sc Benjamin Jourdain} \thanks{Universit\'e Paris-Est, Cermics (ENPC), INRIA, F-77455 Marne-la-Vall\'ee, France. E-mail: {\tt   benjamin.jourdain@enpc.fr}}~\footnotemark[3]
\and   
{\sc  Gilles Pag\`es} \thanks{Laboratoire de Probabilit\'es, Statistique et Mod\'elisation, UMR~8001, Campus Pierre et Marie Curie, Sorbonne Universit\'e case 158, 4, pl. Jussieu, F-75252 Paris Cedex 5, France. E-mail:{\tt  gilles.pages@upmc.fr}}~\thanks{This research
benefited from the support of the ``Chaire Risques Financiers'', Fondation du Risque}  }
\date{}
\maketitle 
\begin{abstract}  Quantization provides a very natural way to preserve the convex order when approximating two ordered probability measures by two finitely supported ones. Indeed, when the convex order dominating original probability measure is compactly supported, it is smaller % in the convex order 
  than any of its dual quantizations while the dominated original measure is greater than any of its stationary (and therefore any of its optimal) quadratic primal quantization. Moreover, the quantization errors then correspond to martingale couplings between each original probability measure and its quantization. This permits to prove that any martingale coupling between the original probability measures can be approximated by a martingale coupling between their quantizations in Wassertein distance with a rate given by the quantization errors but also in the much finer adapted Wassertein distance. As a consequence, while the stability of (Weak) Martingale Optimal Transport problems with respect to the marginal distributions has only been established in dimension $1$ so far, their value function computed numerically for the quantized marginals converges in any dimension to the value for the original probability measures as the numbers of quantization points go to $\infty$.
\end{abstract}

\section*{Introduction}\label{sec:1}

For $d\!\in\N^*$ and $\mu,\nu$ in the set ${\mathcal P}(\R^d)$ of probability measures on $\R^d$, we say that $\mu$ is smaller than $\nu$ in the convex order and denote $\mu\le_{cvx}\nu$ if
\begin{equation}
   \forall \varphi:\R^d\to\R\mbox{ convex },\;\int_{\R^d}\varphi(x)\mu(dx)\le\int_{\R^d}\varphi(y)\nu(dy),\label{eq:defcvx}
\end{equation}
when the integrals make sense (since any real valued convex function is bounded from below by an affine function $\int_{\R^d}\varphi(x)\mu(dx)$ makes sense in $\R\cup\{+\infty\}$ as soon as $\int_{\R^d}|x|\mu(dx)<+\infty$). 
For $p\ge 1$, we denote by ${\mathcal P}_p(\R^d)=\{\mu\!\in{\mathcal P}(\R^d):\int_{\R^d}|x|^p\mu(dx)<+\infty\}$ the Wasserstein space with index $p$ over $\R^d$. When $\mu,\nu\in{\mathcal P}_1(\R^d)$, according to the Strassen theorem~\cite{Strassen},  $\mu\le_{cvx}\nu$ if and only if there exists a martingale coupling between $\mu$ and $\nu$ that is a probability measure $\pi(dx,dy)$ on $\R^d\times\R^d$ with marginals $\int_{y\in\R^d}\pi(dx,dy)$ and $\int_{x\in\R^d}\pi(dx,dy)$  equal to $\mu(dx)$ and $\nu(dy)$ respectively such that $\pi(dx,dy)=\mu(dx)\pi_x(dy)$ for some Markov kernel $\pi_x(dy)$ with the martingale property: $\forall x\!\in\R^d$, $\pi_x\in{\cal P}_1(\R^d)$ and $\int_{\R^d}y\pi_x(dy)=x$. We denote by ${\cal P}(\mu,\nu)$ the set of probability measures on $\R^d\times\R^d$ with respective marginals $\mu$ and $\nu$ and by ${\cal M}(\mu,\nu)$ the subset of ${\cal P}(\mu,\nu)$ consisting of martingale couplings.

Let $(\mu,\nu)$ belong to the set ${\cal P}_\le\times{\cal P}_p(\R^d)$ of couples of elements of ${\cal P}_p(\R^d)$ with the first one smaller than the second in the convex order. In its simplest form, the Martingale Optimal Transport problem consists in computing 
$$V_{c}(\mu,\nu)=\inf_{\pi\in{\cal M}(\mu,\nu)}\int_{\R^d\times\R^d}c(x,y)\pi(dx,dy)$$
and the optimal martingale couplings achieving this infimum for some measurable cost function $c:\R^d\times\R^d\to\R$. When the interest rate is zero, for an exotic option written on $d$ traded assets with payoff given by the function $c$ of their prices at times $s$ and $t$, then $V_c(\mu,\nu)$ (resp. $-V_{-c}(\mu,\nu)$) provides a robust lower (resp. upper) price bound when $\mu$ and $\nu$ are the respective joint laws of these $d$ assets at times $s$ and $t$ (for instance obtained by calibration of a model to vanilla option prices). Since its introduction in~\cite{Beietal}, this MOT problem has received  recently a great attention in the financial mathematics literature. In particular, the structure of martingale optimal transport couplings~\cite{BeJu,CampLaachMart,DeMarchTouzi,Ghousskimlim,HLTo}, continuous time formulations~\cite{DoSo,GalHTo,Toetal}, links with the Skorokhod embedding problem~\cite{Beicoxhues}, numerical methods~\cite{AlJo,AlJo2,DeMarch,GuOb,HL19} and stability properties~\cite{BackPam,JouMar,Wiesel} have been investigated. The MOT problem is a particular instance where the measurable cost function $C:\R^d\times{\cal P}_1(\R^d)\to\R$ is linear in the measure component ($C(x,\eta)=\int_{\R^d}c(x,y)\eta(dy)$) of the Weak Martingale Optimal Transport problem $$\inf_{\pi\in{\cal M}(\mu,\nu)}\int_{\R^d}C(x,\pi_x)\mu(dx)$$
introduced in \cite{BackPam} by adding the martingale constraint to the Weak Optimal Transport problem introduced in \cite{GoRoSamTe}.

To devise a numerical procedure devoted to the computation of the value and of the optimal couplings in the MOT and WMOT problems, a first natural step consists in approximating $\mu$ and $\nu$ by finitely supported probability measures which are still in the convex order. To our  best knowledge, few studies consider the problem of preserving the convex order while approximating a sequence of probability measures. We mention the thesis of Baker~\cite{Baker} who proposes the following construction in dimension $d=1$. Let for $\eta\in{\cal P}(\R)$ and $u \in (0,1)$, $F_{\eta}^{-1}(u)=\inf \{x \in \R: \eta((-\infty,x])\ge u \}$  be the quantile of $\eta$ of order $u$. For $(\mu,\nu)\in{\cal P}_{\le}\times{\cal P}_1(\R)$ and $N,K\in\N^*$ with $N/K\in\N^*$, one has
$$
\frac{1}{N}\sum_{i=1}^{N}\delta_{N\int_{\frac{i-1}{N}}^{\frac{i}{N} } F_\mu^{-1}(u) du }\le_{cvx}\frac{1}{K}\sum_{i=1}^{K}\delta_{K\int_{\frac{i-1}{K}}^{\frac{i}{K} } F_\nu^{-1}(u) du }.
$$
{\em Dual (or Delaunay) quantization} introduced by Pag\`es and Wilbertz ~\cite{PaWi0} and further studied in~\cite{PaWi1,PaWi2, PaWi3} gives another way to preserve the convex order in dimension $d=1$ when using the same grid to quantize both probability measures.

In two recent papers~\cite{AlJo2,AlJo},  Alfonsi, Corbetta and Jourdain propose to restore for $(\mu,\nu)\in{\cal P}_\le\times{\cal P}_1(\R^d)$ the convex ordering from any finitely supported approximations $\tilde \mu$ and $\tilde \nu$ of $\mu$ and $\nu$. In dimension $d=1$, one may define the increasing (resp. decreasing) convex order by adding the constraint that the test function $\varphi$ is non-decreasing (resp. non-increasing) in~\eqref{eq:defcvx}. According to~\cite{AlJo2}, the convex order restauration can be achieved by keeping $\tilde\mu$ (resp. $\tilde\nu$) and replacing $\tilde\nu$ (resp. $\tilde\mu$) by the supremum (resp. infimum) between $\tilde\mu$ and $\tilde\nu$ for the increasing convex order when $\int_\R x\tilde\nu(dx)\le \int_\R x \tilde\mu(dx)$ and the decreasing convex order when $\int_\R x\tilde\nu(dx)\ge \int_\R x \tilde\mu(dx)$. The convex, increasing convex and decreasing convex orders are nicely characterized in terms of the potential function obtained as the anti-derivative of the cumulative distribution function (or of the quantile function). The supremum and infimum of two probability measures for one of these orders can be computed using their potential functions. For a general dimension $d$,~\cite{AlJo} suggests to keep $\tilde \nu$ and replace $\tilde \mu$ by its projection on the set of probability measures dominated by $\tilde \nu$ for the quadratic Wasserstein distance $W_2$ (see \eqref{defwas} below for the definition of this distance). This projection can be computed by solving a quadratic optimization problem with linear constraints.

In the present paper, when $\nu$ is compactly supported, we investigate the combined approximation of $\mu$ by some quadratic-optimal primal quantization and of $\nu$ by some dual quantization. By construction, any quadratic-optimal primal quantization of $\mu$ satisfies a stationarity property which implies that it is smaller than $\mu$ in the convex order. On the other hand, any dual quantization of $\nu$ is greater than this probability measure in the convex order. Therefore the convex order between $\mu$ and $\nu$ is preserved by this combined approximation. Notice that, in contrast with the previous approaches, it cannot be generalized to the convex order preserving approximation of more than two probability measures. Moreover the dual quantization approximation is only possible for compactly supported probability measures. In contrast with these restrictions, we will see that the studied approach proves to provide robust approximations of (Weak) Martingale Optimal Transport problems even in dimension $d\ge 2$.

The first section of the paper is devoted to primal (or Voronoi) quantization. For $\mu\in{\cal P}_p(\R^d)$ with $p\ge 1$, we show that an element of the set ${\cal P}(\R^d,N)$ of probability measures on $\R^d$ whose support contains at most $N$ points is an $L_p$-optimal $N$-quantization of $\mu$ iff it is a $W_p$-projection of $\mu$ on ${\cal P}(\R^d,N)$ where the Wasserstein distance $W_p$ is defined by
\begin{equation}
  W_p(\mu,\nu)^p=\inf_{\pi\in{\cal P}(\mu,\nu)}\int_{\R^d\times\R^d}|x-y|^p\pi(dx,dy).\label{defwas}
\end{equation} In the quadratic $p=2$ case, any stationary and therefore any optimal $N$-quantization $\hat\mu^{N}$ of the measure $\mu$ is smaller than $\mu$ in the convex order. Moreover $W_2(\hat\mu^{N},\mu)=M_2(\hat\mu^{N},\mu)$ where 
\begin{equation}
   M_p(\eta,\nu)=\inf_{\pi\in{\cal M}(\eta,\nu)}\int_{\R^d\times\R^d}|x-y|^p\pi(dx,dy)\mbox{ for }(\eta,\nu)\in{\cal P}_\le\times{\cal P}_1(\R^d).\label{defwasmart}
\end{equation} 
This enables us to check that any quadratic-optimal primal quantizer for $\mu$ remains quadratic-optimal for each probability measure smaller than $\mu$ and greater than its associated quantization for the convex order.

The second section deals with the dual quantization of a compactly supported probability measure $\mu$ which is obtained by minimizing $M_p(\mu,\eta)$ over the subset ${\cal P}_{\ge \mu}(\R^d,N)$ of ${\cal P}(\R^d,N)$ consisting in probability measures greater than $\mu$ in the convex order. Since ${\cal M}(\mu,\eta)\subset{\cal P}(\mu,\eta)$, $M_p(\mu,\eta)\ge W_p(\mu,\eta)$ and, in general, the inequality is strict. It turns out, that, even in the quadratic $p=2$ case, $L_p$-optimal dual quantizations of $\mu$ are not necessarily $W_p$-projections of $\mu$ on ${\cal P}_{\ge \mu}(\R^d,N)$. Nevertheless, we check that any quadratic-optimal dual quantizer of $\mu$ remains quadratic-optimal for each probability measure greater than $\mu$ and smaller than its associated dual quantization in the convex order.

In the third section, we consider two probability measures $\mu,\nu\in{\cal P}(\R^d)$ such that $\mu\le_{cvx}\nu$ with $\nu$ compactly supported. Since the quantization errors between $\mu$ and any of its quadratic-optimal primal $N$-quantization $\hat\mu^N$ and between $\nu$ and any of its $L_p$-optimal dual $K$-quantization $\check\nu^K$  correspond to martingale couplings, we are able to approximate in Wasserstein distance on ${\cal P}(\R^d\times\R^d)$ any martingale coupling $\pi\in{\cal M}(\mu,\nu)$ by a martingale coupling $\bar{\pi}^{N,K}\in{\cal M}(\hat\mu^N,\check\nu^K)$ with a rate given by the quantization errors. We also check that, as $N,K\to\infty$, $\bar{\pi}^{N,K}$ converges to $\pi$ for the much finer adapted Wasserstein distance defined in \eqref{adapwas} below which captures the temporal structure of probability distributions with two time marginals. Numerous financial applications of this adapted Wasserstein distance have been investigated in \cite{BBBE1}. According to \cite{BBBE2}, the topology induced by this distance is equal to the other adapted topologies which had been introduced in particular in view of financial applications.

The adapted Wasserstein is particularly well suited to deal with Weak Martingale Optimal Transport problems. While their stability with respect to the marginal distributions has only been established in dimension $1$ so far \cite{BJMP}, this enables us to check in Section 4 that their value function computed numerically for the quantized marginals converge in any dimension to the value for the original probability measures as the numbers of quantization points go to $\infty$.

For the reader's convenience, we now list the notations and basic properties that have been introduced in the above text.

\medskip
\noindent {\sc Definitions and notations.}

 Let $d\!\in \N^*$

\noindent $\bullet$ $|\cdot|$ denotes the canonical Euclidean norm on $\R^d$.

%\smallskip 

%\smallskip
\noindent $\bullet$ ${\rm conv}(A)$ denotes the (closed) convex hull of $A\subset \R^d$ and ${\rm card}(A)$ or $|A|$  its cardinality (depending on the context).

%\smallskip 
\noindent $\bullet$ Let ${\cal P}(\R^d)$ denote the set of probability measures on $\R^d$ endowed with its Borel sigma-field ${\cal B}(\R^d)$.

%\smallskip 
\noindent $\bullet$ For $p\ge 1$, let ${\cal P}_p(\R^d)=\{\mu\in{\cal P}(\R^d):\int_{\R^d}|x|^p\mu(dx)<\infty\}$.

%\smallskip 
\noindent $\bullet$ For $p\ge 1$, let ${\cal P}_\le\times{\cal P}_p(\R^d)=\{(\mu,\nu)\in{\cal P}_p(\R^d)\times{\cal P}_p(\R^d):\mu\le_{cvx}\nu\}$.

\smallskip
\noindent $\bullet$ For every integer $N\ge 1$, we denote by  $\mathcal{P}(\R^d, N)$ the set of distributions on $\R^d$ whose support contains at most $N$ points.

%\smallskip
\noindent $\bullet$ For $\mu,\nu\in{\cal P}(\R^d)$, let $${\cal P}(\mu,\nu)=\left\{\pi\in{\cal P}(\R^{d}\times\R^{d}):\;\forall A\in{\cal B}(\R^d),\;\pi(A\times\R^d)=\mu(A)\mbox{ and }\pi(\R^d\times A)=\nu(A)\right\}$$ denote the set of couplings between $\mu$ and $\nu$. For $\pi\in{\cal P}(\mu,\nu)$, we denote by $\pi_x(dy)$ the ($\mu(dx)$ a.e. unique) Markov kernel such that $\pi(dx,dy)=\mu(dx)\pi_x(dy)$.

%\smallskip
\noindent $\bullet$ For $\mu,\nu\in{\cal P}_1(\R^d)$ with $\mu\le_{cvx}\nu$, let
$${\cal M}(\mu,\nu)=\left\{\pi\in{\cal P}(\mu,\nu):\;\mu(dx)\mbox{ a.e. }\int_{\R^d}y\pi_x(dy)=x\right\}$$
denote the (non empty according to the Strassen theorem) set of martingale couplings between $\mu$ and $\nu$. 

%\smallskip
\noindent $\bullet$  For $p\ge 1$ and $\mu,\nu\in{\cal P}(\R^d)$, let \begin{equation*}
   \cW_p(\mu, \nu)=\inf_{\pi\in{\cal P}(\mu,\nu)}\left(\int_{\R^d\times\R^d}|x-y|^p \pi(dx,dy)\right)^{1/p}\le\infty 
\end{equation*} denote the Wasserstein distance with index $p$. This is a complete metric on $\mathcal{P}_p(\R^d)$.

%\smallskip
\noindent $\bullet$  For $p\ge 1$ and $(\mu,\nu)\in{\cal P}_\le\times{\cal P}_1(\R^d)$, let \begin{equation*}
   \cM_p(\mu, \nu)=\inf_{\pi\in{\cal M}(\mu,\nu)}\left(\int_{\R^d\times\R^d}|x-y|^p \pi(dx,dy)\right)^{1/p}\le\infty.
\end{equation*}
Since ${\cal M}(\mu,\nu)\subset{\cal P}(\mu,\nu)$, we clearly have, $\cW_p(\mu,\nu)\le\cM_p(\mu,\nu)$. In the quadratic $p=2$ case, when $\mu,\nu\in{\cal P}_2(\R^d)$,
\begin{align*}
  \forall\pi\in{\cal M}(\mu,\nu),\;\int_{\R^d\times\R^d}|y-x|^2\pi(dx,dy)&=\int_{\R^d}|y|^2\nu(dy)-2\int_{\R^d}x.\int_{\R^d}y\pi_x(dy)\mu(dx)+\int_{\R^d}|x|^2\mu(dx)\\&=\int_{\R^d}|y|^2\nu(dy)-\int_{\R^d}|x|^2\mu(dx)
\end{align*} so that \begin{equation}
   M_2^2(\mu,\nu)=\int_{\R^d}|y|^2\nu(dy)-\int_{\R^d}|x|^2\mu(dx).\label{m2mom}
\end{equation}
%\smallskip

\noindent $\bullet$ For $\pi\in{\cal P}(\mu,\nu)$ and $\tilde\pi\in{\cal P}(\tilde\mu,\tilde\nu)$ we consider the adapted Wasserstein distance with index $p\ge 1$ between $\pi(dx,dy)=\mu(dx)\pi_x(dy)$ and $\tilde\pi(d\tilde x,d\tilde y)=\tilde\mu(d\tilde x)\tilde\pi_{\tilde x}(d\tilde y)$ :
   \begin{align}
\cAW_p(\pi,\tilde \pi)=\inf_{m\in{\cal P}(\mu,\tilde\mu)}\left(\int_{\R^d\times\R^d}(|x-\tilde x|^p+\cW_p^p(\pi_x,\tilde\pi_{\tilde x}))m(dx,d\tilde x)\right)^{1/p}.\label{adapwas}
   \end{align}

\section{Primal (Voronoi) quantization and Wasserstein projection} 
In this subsection, we make a connection between primal quantization and various projections (in the Wasserstein sense), including, in the quadratic case, with the one mentioned above in the introduction.
Let us first recall the following basic  facts about the (primal) Voronoi quantization of $\mu\in{\cal P}_p(\R^d)$ with $p\ge 1$ (see~\cite{GrLu, Pag2015, PagSpring2018} among others):

\smallskip
-- Let $\Gamma= \{x_1, \ldots,x_{_N}\}\subset \R^d$ denote a finite subset of size $N$. The $L^p$-quantization error modulus  $e_{p}(\Gamma, \mu)$ satisfies 
\begin{equation}\label{eq:2.1.0}
e_{p}(\Gamma, \mu)^p  = \int_{\R^d} |x-{\rm Proj}_{\Gamma}(x)|^p\mu(dx)
\end{equation}
where ${\rm Proj}_{\Gamma}$ denotes a  Borel  nearest neighbour projection on $\Gamma$ satisfying $|x-{\rm Proj}_{\Gamma}(x)|={\rm dist}(x,\Gamma)$. If $X\sim\mu$, the random variable ${\rm Proj}_{\Gamma}(X)$ with law $\mu\circ{\rm Proj}_{\Gamma}^{-1}$ is called $\Gamma$-quantization of $X$.

\smallskip
--  For any {\em level} $N\ge 1$, there exists an optimal {\em grid} or {\em $N$-quantizer} $\Gamma_{p,N}$ such that 
\[
e_{p,N}(\mu):= \inf\Big\{e_p(\Gamma,\mu): \Gamma \subset \R^d, \; |\Gamma| \le N\Big\}= e_p\big(\Gamma_{p,N},\mu\big).
\]
When ${\rm card}({\rm supp}(\mu))\le N$ then $\Gamma_{p,N}={\rm supp}(\mu)$ and when ${\rm card}({\rm supp}(\mu))> N$, then $\Gamma_{p,N}$ has exactly $N$ pairwise distinct elements. Moreover, Theorem~4.2 in~\cite{GrLu} ensures that \begin{equation}
   \mu\Big(\{x\in\R^d:\exists y\neq \tilde y\in \Gamma_{p,N},\;|x-y|=|x-\tilde y|={\rm dist}(x,\Gamma_{p,N})\}\Big)=0,\label{munechargepasfront}
 \end{equation} so that $\mu\circ{\rm Proj}_{\Gamma_{p,N}}^{-1}$ does not depend on the choice of the Borel nearest neighbour projection ${\rm Proj}_{\Gamma_{p,N}}$ on $\Gamma_{p,N}$. We denote this probability measure by $\hat\mu^{\Gamma_{p,N}}$. In the same way, when $X\sim\mu$, ${\rm Proj}_{\Gamma_{p,N}}(X)$ a.s. does not depend on the Borel nearest neighbour projection ${\rm Proj}_{\Gamma_{p,N}}$ and is denoted $\hat X^{\Gamma_{p,N}}$.

\smallskip
-- In the quadratic case ($p=2$), any optimal quantization grid $\Gamma_{2,N}$ (possibly not unique) and its induced quantization $\hat X^{\Gamma_{2,N}}= {\rm Proj}_{\Gamma_{2,N}}(X)$ with distribution $\hat \mu^{\Gamma_{2,N}}$ satisfy a {\em stationarity} (or {\em self-consistency}) property (see e.g.~\cite{GrLu},~\cite{Pag2015} or~\cite{PagSpring2018}, Proposition~5.1 among others) that is
\begin{equation}\label{eq:VoroStatio}
\E\big(X\,|\, \hat X^{\Gamma_{2,N}}\big) =  \hat X^{\Gamma_{2,N}}\mbox{ so that }\hat \mu ^{\Gamma_{2,N}} \le_{cvx} \mu.
\end{equation}
Indeed, the support of the distribution of $\E\big(X\,|\, \hat X^{\Gamma_{2,N}}\big)$ is equal to $\tilde \Gamma_{2,N}:=\{\E(X|\hat X^{\Gamma_{2,N}}=x),x\in \Gamma_{2,N}\}$ (when ${\rm card}({\rm supp}(\mu))\le N$, $\hat X^{\Gamma_{2,N}}=X$ and $\tilde \Gamma_{2,N}=\Gamma_{2,N}={\rm supp}(\mu)$), contains at most $N$ points and $e_2(\tilde \Gamma_{2,N},\mu)^2=\E({\rm dist}(X,\tilde\Gamma_{2,N})^2)\le \E(|X-\E\big(X\,|\, \hat X^{\Gamma_{2,N}}\big)|^2)$. As a consequence,
\begin{align*}
  e_{2,N}(\mu)^2&\le e_2(\tilde \Gamma_{2,N},\mu)^2\le \E(|X-\E\big(X\,|\, \hat X^{\Gamma_{2,N}}\big)|^2)=\E(|X-\hat X^{\Gamma_{2,N}}|^2)
-\E(|\hat X^{\Gamma_{2,N}}-\E\big(X\,|\, \hat X^{\Gamma_{2,N}}\big)|^2)\\&=e_{2,N}(\mu)^2-\E(|\hat X^{\Gamma_{2,N}}-\E\big(X\,|\, \hat X^{\Gamma_{2,N}}\big)|^2),
\end{align*}
so that the second term in the right-hand side vanishes. Therefore the distribution of $(\hat{X}^{\Gamma_{2,N}},X)$ belongs to ${\cal M}(\hat\mu^{\Gamma_{2,N}},\mu)$ and
\begin{equation}
   e_{2,N}^2(\mu)=\E[|X-\hat{X}^{\Gamma_{2,N}}|^2]=\E[|X|^2]-\E[|\hat{X}^{\Gamma_{2,N}}|^2]=M_2^2(\hat\mu^{\Gamma_{2,N}},\mu).\label{e2m2}
\end{equation}

% We start by a lemma which makes a first connection between Wasserstein distance and quantization error.

\begin{Proposition}Let $p\!\in [1, +\infty)$ and $\mu \!\in {\cal P}_p(\R^d)$.\label{propprimcvx}

\smallskip 
\noindent    $(a)$  Let $\Gamma\subset \R^d$ be a finite set and ${\cal P}(\Gamma)$ denote the subset of $\Gamma$-supported distributions. Then
\[
\cW_p\big(\mu,{\cal P}(\Gamma)\big) := \inf_{\nu \in{\cal P}(\Gamma)}\cW_p(\mu,\nu)  =  e_{p}(\Gamma, \mu):= \big\|{\rm dist}(., \Gamma)\big\|_{L^p(\mu)}
\]
and for any Borel  nearest neighbour projection ${\rm Proj}_{\Gamma}$ on $\Gamma$, $\mu \circ {\rm Proj}_{\Gamma}^{-1}$ is a $\cW_p$-projection of $\mu$ on ${\cal P}(\Gamma)$.

\smallskip 
\noindent $(b)$ The probability measure $\nu\in{\cal P}(\R^d,N)$ is a $W_p$-projection of $\mu$ on ${\cal P}(\R^d,N)$ iff $\nu=\hat\mu^{\Gamma_{N}}$ for some $L_p$-optimal $N$-quantizer $\Gamma_{N}$  of $\mu$. Moreover, $W_p(\mu,{\cal P}(\R^d,N))=e_{p,N}(\mu)$.

\smallskip 
\noindent $(c)$ Quadratic case ($p=2$). A subset $\Gamma$ of $\R^d$ with cardinality at most $N$ is a quadratic optimal $N$-quantizer of $\mu$ iff there exists a probability measure $\nu\in{\cal P}_{\le\mu} (\R^d,N)$ such that $\nu(\Gamma)=1$ and one of the following equivalent conditions is satisfied
\begin{itemize}
\item $\nu$ is a $\cW_2$-projection of $\mu$ on ${\cal P}_{\le\mu} (\R^d,N)$ i.e. $$\cW_2(\mu,\nu)=\cW_2(\mu,{\cal P}_{\le\mu} (\R^d,N)):=\inf_{\eta\in{\cal P}_{\le\mu} (\R^d,N)}\cW_2(\mu,\eta),$$
  \item $\int_{\R^d}|x|^2\nu(dx)=\sup_{\eta\in{\cal P}_{\le\mu} (\R^d,N)}\int_{\R^d}|x|^2\eta(dx)$.
  \end{itemize}
  Moreover, we then have $W_2(\nu,\mu)=M_2(\nu,\mu)$ and $\nu=\hat\mu^\Gamma$.
\end{Proposition}
Apart from the interpretation in terms of $\cW_p$-projection, the first statement can be found in Lemma~3.4 p.33 \cite{GrLu}. 
Before proving the proposition, let us state and check some easy consequence of the necessary and sufficient condition in $(b)$. 
\begin{Corollary}\label{coroptgridprim}
   Let $\Gamma_{2,N}$ be a quadratic optimal $N$-quantizer of $\mu\in{\cal P}_2(\R^d)$. Then for any probability measure $\nu$ such that $\hat{\mu}^{\Gamma_{2,N}}\le_{cvx}\nu\le_{cvx}\mu$, $\Gamma_{2,N}$ is a quadratic optimal $N$-quantizer of $\nu$ and $\hat{\nu}^{\Gamma_{2,N}}=\hat{\mu}^{\Gamma_{2,N}}$.
 \end{Corollary} \noindent {\bf Proof of Corollary~\ref{coroptgridprim}.} Since $\hat{\mu}^{\Gamma_{2,N}}\in {\cal P}_{\le\nu} (\R^d,N)\subset{\cal P}_{\le\mu} (\R^d,N)$ and $\int_{\R^d}|x|^2\hat{\mu}^{\Gamma_{2,N}}(dx)=\sup_{\eta\in{\cal P}_{\le\mu} (\R^d,N)}\int_{\R^d}|x|^2\eta(dx)$ by the necessary condition in Proposition~\ref{propprimcvx} $(c)$, one has $$\int_{\R^d}|x|^2\hat{\mu}^{\Gamma_{2,N}}(dx)=\sup_{\eta\in{\cal P}_{\le\nu} (\R^d,N)}\int_{\R^d}|x|^2\eta(dx).$$  Therefore, by the sufficient condition in Proposition~\ref{propprimcvx} $(c)$, $\Gamma_{2,N}$ is a quadratic optimal $N$-quantizer of $\nu$ and $\hat{\nu}^{\Gamma_{2,N}}=\hat{\mu}^{\Gamma_{2,N}}$.
 \hfill$\Box$

\bigskip

\noindent {\bf Proof of Proposition~\ref{propprimcvx}.} $(a)$ Let $\nu\!\in {\cal P}(\Gamma)$ and $\pi\in{\cal P}(\mu,\nu)$. Then $\mu(dx)$ a.e., $\pi_x(\Gamma)=1$ so that $\pi_x(dy)$ a.e. $|x-y|\ge{\rm dist}(x,\Gamma)$. Therefore
\begin{align*}
\int_{\R^d\times\R^d} |x-y|^p \pi(dx,dy) &\ge \int{\rm dist}(x, \Gamma)^p \mu(dx) = e_p(\Gamma, \mu)^p.
\end{align*}
Taking the infimum over $\pi\in {\cal P}(\mu,\nu)$ and $\nu\!\in {\cal P}(\Gamma)$, we deduce that $\cW_p\big(\mu,{\cal P}(\Gamma)\big)^p\ge e_p(\Gamma, \mu)^p$.
Now let ${\rm Proj}_{\Gamma}$ denote a Borel nearest neighbour projection on $\Gamma$. Since $\mu \circ {\rm Proj}_{\Gamma}^{-1}\in{\cal P}(\Gamma)$ and $\mu\circ(I_d,{\rm Proj}_{\Gamma})^{-1}\in{\cal P}(\mu,\mu \circ {\rm Proj}_{\Gamma}^{-1})$, we have
$$
\cW^p_p\big(\mu,{\cal P}(\Gamma)\big) \le \cW^p_p(\mu, \mu \circ {\rm Proj}_{\Gamma}^{-1})\le  \int_{\R^d} |x-{\rm Proj}_{\Gamma}(x)|^p\mu(dx)=e_{p}(\Gamma, \mu)^p,
$$
where the equality follows from~\eqref{eq:2.1.0}. 
Therefore the inequalities are equalities and $\mu \circ {\rm Proj}_{\Gamma}^{-1}$ is a $\cW_p$-projection of $\mu$ on ${\cal P}(\Gamma)$.
 
\smallskip
\noindent $(b)$ Let $\Gamma_{p,N}$ be an $L^p$-optimal $N$-quantizer of $\mu$. Then, for any subset $\Gamma$ of $\R^d$ with at most $N$ points, $e_p(\Gamma_{p,N},\mu)=e_{p,N}(\mu)\le e_p(\Gamma,\mu)$. With $(a)$, we deduce that
$$W_p(\mu,\hat\mu^{\Gamma_{p,N}})\le W_p(\mu,{\cal P}(\Gamma)).$$ By taking the infimum over $\Gamma$, we deduce that $W_p(\mu,\hat\mu^{\Gamma_{p,N}})\le W_p(\mu,{\cal P}(\R^d,N))$ and $\hat\mu^{\Gamma_{p,N}}$ is a $W_p$-projection of $\mu$ on ${\cal P}(\R^d,N)$ so that
$$e_{p,N}(\mu)=e_p(\mu,\Gamma_{p,N})=W_p(\mu,\hat\mu^{\Gamma_{p,N}})=W_p(\mu,{\cal P}(\R^d,N)).$$

Conversely, let $\nu$ be a projection of $\mu$ on ${\cal P}(\R^d,N)$ 
and let $\Gamma_N=\{x\in\R^d:\nu(\{x\})>0\}$. The cardinality of $\Gamma_N$ is at most $N$. For any subset $\Gamma$ of $\R^d$ with at most $N$ points, \begin{equation}
   W_p(\mu,{\cal P}(\Gamma_N))\le W_p(\mu,\nu)=W_p(\mu,{\cal P}(\R^d,N))\le W_p(\mu,{\cal P}(\Gamma)).\label{projpgamn}
\end{equation}
Therefore, by $(a)$, $e_p(\Gamma_N,\mu)\le e_p(\Gamma,\mu)$ and since $\Gamma$ is arbitrary, we deduce that $\Gamma_N$ is an $L^p$-optimal $N$-quantizer of $\mu$. Moreover, the choice $\Gamma=\Gamma_N$ in \eqref{projpgamn} implies that the first inequality is an equality so that, with $(a)$, $W^p_p(\mu,\nu)=\int_{\R^d}{\rm dist}(x,\Gamma_N)^p\mu(dx)$. Hence, for any $W_p$-optimal coupling $\pi\in{\cal P}(\mu,\nu)$, $$\int_{\R^d\times\R^d}(|y-x|^p-{\rm dist}(x,\Gamma_N)^p)\pi(dx,dy)=0.$$
Since $1=\nu(\Gamma_N)=\int_{\R^d}\pi_x(\Gamma_N)\mu(dx)$, $\pi(dx,dy)$ a.e., $|y-x|^p\ge {\rm dist}(x,\Gamma_N)^p$. Therefore, $\pi(dx,dy)$ a.e. $y\in\Gamma_N$ and $|y-x|={\rm dist}(x,\Gamma_N)$. With \eqref{munechargepasfront}, we conclude that
$\mu(dx)$ a.e. there is a unique point $\hat x^{\Gamma_N}\in\Gamma_N$ such that $|\hat x^{\Gamma_N}-x|= {\rm dist}(x,\Gamma_N)$ and $\pi_x(dy)=\delta_{\hat x^{\Gamma_N}}(dy)$. Therefore the second marginal $\nu$ of $\pi$ is equal to $\hat\mu^{\Gamma_N}$.

\smallskip
\noindent $(c)$ In the quadratic case ($p=2$), by $(b)$, \eqref{eq:VoroStatio} and \eqref{e2m2}, any $W_2$-projection $\nu$ of $\mu$ on ${\cal P}(\R^d,N)$ is characterized by the existence of a quadratic optimal $N$-quantizer $\Gamma_{2,N}$ of $\mu$ such that $\nu=\hat\mu^{\Gamma_{2,N}}$, belongs to the smaller set ${\cal P}_{\le \mu}(\R^d,N)$ and satisfies $e_{2,N}^2(\mu)=M_2^2(\nu,\mu)$. Therefore the $W_2$-projections of $\mu$ on ${\cal P}(\R^d,N)$ and on ${\cal P}_{\le \mu}(\R^d,N)$ coincide. To conclude the proof, let us check that $\nu$ is such a projection iff $\int_{\R^d}|x|^2\nu(dx)=\sup_{\eta\in{\cal P}_{\le\mu} (\R^d,N)}\int_{\R^d}|x|^2\eta(dx)$. 
If $\eta\in{\cal P}_{\le\mu} (\R^d,N)$, then, by the comparison between $W_2$ and $M_2$  given in the introduction and \eqref{m2mom}, 
one has\begin{equation}
   \cW_2^2(\mu,\eta)\le\cM_2^2(\eta,\mu)=\int_{\R^d}|y|^2\mu(dy)-\int_{\R^d}|x|^2\eta(dx).\label{ineqw2mart}
\end{equation}
Let $\nu$ be a $W_2$-projection of $\mu$ on ${\cal P}_{\le \mu}(\R^d,N)$. Using $(b)$ for the third equality then \eqref{ineqw2mart} for the last inequality and the last equality, we obtain that
 \begin{align*}
 \int_{\R^d}|y|^2\mu(dy)-\int_{\R^d}|x|^2\nu(dx)&=M_2^2(\nu,\mu)= e_{2,N} (\mu)^2=\cW^2_2\big(\mu,{\cal P}(\R^d, N) \big) \\&=\cW^2_2\big(\mu,{\cal P}_{\le\mu}(\R^d, N) \big)\le\inf_{\eta\in {\cal P}_{\le\mu}(\R^d, N)}\cM^2_2(\eta,\mu)\\&=\int_{\R^d}|y|^2\mu(dy)-\sup_{\eta\in{\cal P}_{\le\mu} (\R^d,N)}\int_{\R^d}|x|^2\eta(dx).
 \end{align*}
 Since $\nu \!\in {\cal P}_{\le\mu}(\R^d,N)$, the two inequalities are equalities. Therefore $$\int_{\R^d}|x|^2\nu(dx)=\sup_{\eta\in{\cal P}_{\le\mu} (\R^d,N)}\int_{\R^d}|x|^2\eta(dx)$$ and since $W^2_2(\mu, {\cal P}(\R^d, N) \big)\le W^2_2(\mu,\nu)\le M_2^2(\nu,\mu)$, these two inequalities are equalities and $W^2_2(\mu,\nu)=M_2^2(\nu,\mu)$. Moreover,
 \begin{equation*}
   e_{2,N} (\mu)^2=\cW^2_2\big(\mu,{\cal P}_{\le\mu}(\R^d, N) \big)=\int_{\R^d}|y|^2\mu(dy)-\sup_{\eta\in{\cal P}_{\le\mu} (\R^d,N)}\int_{\R^d}|x|^2\eta(dx).
 \end{equation*}
If $\nu \in {\cal P}_{\le\mu}(\R^d,N)$ is such that $\int_{\R^d}|x|^2\nu(dx)=\sup_{\eta\in{\cal P}_{\le\mu} (\R^d,N)}\int_{\R^d}|x|^2\eta(dx)$, the last equality combined with \eqref{ineqw2mart} written for $\eta=\nu$ ensures that  $\nu$ is a $\cW_2$-projection of $\mu$ on ${\cal P}_{\le\mu}(\R^d,N)$. \hfill$\Box$

\bigskip
\noindent {\bf Remark about uniqueness.} As a consequence of Proposition~\ref{propprimcvx} $(b)$, it turns out that the uniqueness of $\cW_p$-projections of $\mu$ on ${\cal P}(\R^d,N)$, that of distributions $\hat \mu^N$ of $L^p$-optimal $N$-quantizations and that of $L^p$-optimal $N$-quantizers are equivalent. In dimension $d=1$, for $p=2$, distributions  with $\log$-concave densities have a unique optimal $N$-quantizer (see Kiefer~\cite{Kieff}) hence this projection is unique. In higher dimension,  a general result seems difficult to reach: indeed, the ${\cal N}(0;I_d)$ distribution, being invariant under the action of ${\cal O}(d, \R)$ (orthogonal transforms), so are the (hence infinite)  sets of  its optimal quantizers at levels $N\ge 2$. 
% \textcolor{blue}{Benjamin, how what we know on $\cW_p$-projection can improve this microscopic result?}

Let us recall the sharp rate of convergence of the $L^p$-quantization error stated for instance in Theorem 5.2 \cite{PagSpring2018}.
\begin{Theorem}[Pierce Lemma for primal quantization]\label{thm:zadoretpierce1} 
% \noindent $(a)$ \noindent {\sc Zador's Theorem  for (primal) Voronoi quantization:} Let $X\!\in L_{\R^d}^{p+\eta}(\Omega,{\cal A}, \P)$, $p,\eta>0$,  be a random vector with distribution $\P_{_{\!X}}= \varphi.\lambda_d\stackrel{\perp}{+}\nu_{_{\!X}}$ where $\lambda_d$ denotes the Lebesgue measure and $\nu_{_X}$  denotes the singular part of the distribution. Then
% \[
% \lim_{N\to+\infty} N^{\frac 1d} e_{p,N}(X) = \widetilde J^{vor}_{d,p}\left(\int_{\R^d} \varphi^{\frac{d}{d+p}}d\lambda_d\right)^{\frac 1d +\frac 1p}
% \]
% where $\widetilde J^{vor}_{d,p}=\inf_{N\ge 1} N^{\frac 1d} d_{p,N}\big(\mathcal{U}([0,1]^d)\big)$.
%($J_{p,d}$ is the similar  constant from the original Zador's Theorem). \textcolor{red}{ATTENTION! Les notations $J_{p,d}$ c'est pour la puissance p ieme, je crois\dots)}
%
%\smallskip
%When $d=1$, $\widetilde J^{vor}_{1,p}=\frac{1}{2(p+1)^{1/p}}$.

% \medskip
% \noindent $(b)$ {\sc Non-asympotic  bound (Pierce lemma):} 
  Let $p\ge 1$ and $\eta >0$. For every dimension $d\ge 1$, there exists a real constant $\widetilde C^{vor}_{d,\eta,p} >0$ such that, for every random vector $X:(\Omega,{\cal A}, \P) \to \R^d$,
  \begin{equation}
  e_{p,N}(X)\le \widetilde C^{vor}_{d,\eta,p} N^{-\frac 1d} \sigma_{p+\eta}(X) \label{majopierceprim}
  \end{equation}
where, for every $r>0$, $\sigma_r(X)= \inf_{a\in \R^d}\|X-a\|_r\le +\infty$.
\end{Theorem}

\noindent {\bf Example.} Let us explicit the quadratic optimal primal $N$-quantizer of $\mu(dx)=\frac{1_{]0,1[}(x)}{2\sqrt{x}}dx$. In \cite{FortPag}, by checking that the distortion function has a unique critical point, Fort and Pagès prove uniqueness and derive semi-closed forms for optimal quantizers of three families of one-dimensional distributions indexed by $\rho>0$ : the exponential distributions $1_{(0,\infty)}(x)\rho e^{-\rho x}dx$, the power distributions $\rho 1_{(0,1)}(x)x^{\rho-1}dx$  on the interval $(0,1)$ which include $\mu$ for $\rho=\frac 12$ and the power distributions $\rho 1_{(1,\infty)}(x)x^{-(1+\rho)}dx$ on the interval $(1,+\infty)$. For $\mu(dx)$ and $p=2$, their semi-closed form relies on the inductive solution of third degree equations. Working with a different parametrization, we will end up with much simpler second degree equations. We have $F_\mu(x)=1_{\{x\ge 0\}}\sqrt{x}$ and $F_\mu^{-1}(u)=u^2$ for $u\in (0,1)$. We parametrize the boundaries of the optimal quadratic Voronoi cells by $x_{k+\frac{1}{2}}=F_\mu^{-1}(q_k)=q_k^2$ for $k\in\{0,\hdots,N\}$ where $0=q_0<q_1<q_2<\hdots<q_{N-1}<q_N=1$. In particular $x_{0+\frac12}=0$ and $x_{N+\frac12}=1$. By the stationarity property \eqref{eq:VoroStatio}, for $k\in\{1,\hdots,N\}$, $x_k=\frac{1}{q_k-q_{k-1}}\int_{x_{k-\frac{1}{2}}}^{x_{k+\frac{1}{2}}}x\mu(dx)$ (by \eqref{munechargepasfront}, $\mu(\{x_{k-\frac{1}{2}}\})=\mu(\{x_{k+\frac{1}{2}}\})=0$ and we do not need to worry about the inclusion of the boundary points in the integration interval). The equality $x_{k+\frac{1}{2}}=\frac{x_k+x_{k+1}}{2}$ valid for $k\in\{1,\hdots, N-1\}$ also writes
$$F_\mu^{-1}(q_k)=\frac{1}{2}\left(\frac{1}{q_k-q_{k-1}}\int_{q_{k-1}}^{q_k}F_\mu^{-1}(u)du+\frac{1}{q_{k+1}-q_{k}}\int_{q_{k}}^{q_{k+1}}F_\mu^{-1}(u)du\right).$$
For the above choice of $\mu$, we obtain
$$q_k^2=\frac{1}{6}\left(q_k^2+q_kq_{k-1}+q_{k-1}^2+q_{k+1}^2+q_{k+1}q_k+q_{k}^2
\right)\mbox{ so that }q_{k+1}^2+q_{k+1}q_k+q_kq_{k-1}+q_{k-1}^2-4q_k^2=0.$$
We deduce that $q_{k+1}=c_{k+1}q_1$, where $c_{k+1}^2+c_{k+1}c_k+c_kc_{k-1}+c_{k-1}^2-4c_k^2=0$ so that $c_{k+1}=\frac{\sqrt{17 c_k^2-4c_kc_{k-1}-4c_{k-1}^2}-c_k}{2}$. Starting from $c_0=0$ and $c_1=1$, we easily compute inductively the factors $c_k$ (for instance $c_2=\frac{\sqrt{17}-1}{2}$). To ensure $q_N=1$, we need $q_1=\frac{1}{c_N}$. Therefore $q_k=\frac{c_k}{c_N}$ for $k\in\{0,\hdots,N\}$ and the optimal quadratic primal grid is given by
$$x_k=\frac{1}{q_k-q_{k-1}}\int_{q_{k-1}}^{q_k}F_\mu^{-1}(u)du=\frac{q_k^2+q_kq_{k-1}+q_{k-1}^2}{3}=\frac{c_k^2+c_kc_{k-1}+c_{k-1}^2}{3c_N^2},\;k\in\{1,\hdots,N\}.$$
Of course for $a<b$, when $\mu(dx)=\frac{1_{]a,b[}(x)}{2\sqrt{(b-a)(x-a)}}dx$, then $x_k=a+(b-a)\frac{c_k^2+c_kc_{k-1}+c_{k-1}^2}{3c_N^2}$ and when $\mu(dx)=\frac{1_{]a,b[}(x)}{2\sqrt{(b-a)(b-x)}}dx$, $x_{N+1-k}=b-(b-a)\frac{c_k^2+c_kc_{k-1}+c_{k-1}^2}{3c_N^2}$.
\section{Dual (Delaunay) quantization} \label{subsec:dualQ}
We assume throughout this section that $\mu$ is compactly supported. Let $X:(\Omega, {\cal A}, \P)\to \R^d$ be a random vector lying in $L^{\infty}(\P)$ with distribution $\mu$. Optimal dual (or Delaunay) quantization as introduced in~\cite{PaWi1} relies on the best
approximation which can be achieved by a discrete random vector $\check X$ that
satisfies a certain stationarity assumption on the extended probability space
$(\Omega\times [0,1], {\cal A}\otimes {\cal B}([0,1]), \P\otimes \lambda)$  where ${\cal B}([0,1])$ and $\lambda$ respectively denote the Borel sigma field and the Lebesgue measure on the interval $[0,1]$. 
To be more precise,  we define, for $p\!\in [1,+\infty)$,  
\begin{eqnarray*}
d_{p,N}(X) & = &  \inf_{ \check X}\Big\{ \big\|X -  \check X\big\|_p: \check X:(\Omega\times [0,1], {\cal A}\otimes {\cal B}([0,1]), \P\otimes \lambda)\to
\R^d,  \\
& & \qquad\qquad\qquad\qquad{\rm card}\big(\check
X(\Omega\times[0,1])\big) \leq N \text{ and } \E(\check X|X) = X \Big\}.
\end{eqnarray*}
For every level $N\ge d+1$, the set of such $\check X$ is not empty. Indeed, one may choose $d+1$ points whose convex hull has a non empty interior and includes the support of $\mu$. Then the unique probability measure supported on these points with the same expectation as $\mu$ belongs to the set ${\mathcal P}_{\ge \mu}(\R^d,d+1)$ of distributions dominating $\mu$ for the convex order and supported by at  most $d+1$ elements. By Lemma 2.22 in~\cite{Kallenberg}, we see that for each $\nu\in{\mathcal P}_{\ge \mu}(\R^d,N)$ and each martingale coupling $\pi\in{\cal M}(\mu,\nu)$, there exists on $(\Omega\times [0,1], {\cal A}\otimes {\cal B}([0,1]), \P\otimes \lambda)$ a random vector $\check{X}$ such that $(X,\check{X})$ is distributed according to $\pi$ and therefore satisfies $\E(\check X|X) = X$. Hence
\begin{equation}
   d_{p,N}(X)^p=\inf_{\nu\in{\mathcal P}_{\ge \mu}(\R^d,N)}\inf_{\pi\in{\mathcal M}(\mu,\nu)}\int_{\R^d\times\R^d}|y-x|^p \pi(dx,dy)=\inf_{\nu\in{\mathcal P}_{\ge \mu}(\R^d,N)}\cM_p^p(\mu,\nu).\label{dpN}
\end{equation}
As a consequence, $d_{p,N}(X)$ only depends on the distribution $\mu$ of $X$ and can subsequently also be denoted $d_{p,N}(\mu)$.
Next, one easily checks that ${\mathcal P}_{\ge \mu}(\R^d,N)=\bigcup_{\Gamma\in{\cal G}_N}{\mathcal P}_{\ge \mu}(\Gamma)$ where
\begin{align*}
  &{\cal G}_N=\{\Gamma\subset\R^d\mbox{ with cardinality $\le N$ and such that } {\rm supp}(\mu) \subset {\rm conv}(\Gamma)\},\\
  &\mbox{ and }{\mathcal P}_{\ge \mu}(\Gamma)=\{\nu\in{\cal P}(\R^d):\mu\le_{cvx}\nu\mbox{ and }\nu(\Gamma)=1\}.
\end{align*}
For $\Gamma\in{\cal G}_N$, there exists a dual projection
${\rm Proj}_{\Gamma}^{del}: {\rm conv}\big(\Gamma\big)\times [0,1]\to \Gamma$, also called  a {\em splitting operator}, which satisfies, beyond measurability, the following stationarity property
\begin{equation}\label{eq:Statio1}
\forall\, y\!\in {\rm conv}\big(\Gamma\big), \quad \int_0^1{\rm Proj}_{\Gamma}^{del}(y,u) du = y,
\end{equation}
from which one derives the  dual stationarity property 
 \begin{equation}\label{eq:StatioDual1}
 \E\,\big({\rm Proj}^{del}_{ \Gamma}(X,U)\,\big|\,X\big) = X\mbox{ when $U\sim \mathcal{U}([0,1])$ is independent of $X$}. 
\end{equation}
The stationarity property remains valid as soon as $X$ is ${\rm conv}(\Gamma)$-valued and implies that the distribution of ${\rm Proj}^{del}_{ \Gamma}(X,U)$ belongs to ${\mathcal P}_{\ge \mu}(\Gamma)$ which is therefore non empty.

For $\Gamma\in{\cal G}_N$, let $$d_p(\mu,\Gamma)^p=\inf_{\nu\in{\mathcal P}_{\ge \mu}(\Gamma)}\cM_p^p(\mu,\nu),$$
so that $d_{p,N}(\mu)^p=\inf_{\Gamma\in{\cal G}_N}d_p(\mu,\Gamma)$.

In dimension $d=1$, when $\Gamma=\{x_1,x_2,\hdots,x_N\}$ with $x_1<x_2<\hdots<x_N$, the probability measure minimizing $\nu\mapsto\cM_p^p(\mu,\nu)$ over ${\mathcal P}_{\ge \mu}(\Gamma)$ is the distribution $\check \mu^\Gamma$ of ${\rm Proj}_{\Gamma}^{del}(X,U)$ for the splitting operator
$${\rm Proj}_{\Gamma}^{del}(x,u)=\sum_{i=1}^{N-1}1_{[x_i,x_{i+1})}(x)\left(1_{\{u\le \frac{x_{i+1}-x}{x_{i+1}-x_i}\}}x_i+1_{\{u>\frac{x_{i+1}-x}{x_{i+1}-x_i}\}}x_{i+1}\right)+1_{\{x=x_N\}}x_N.$$
Moreover, the coupling minimizing $\int_{\R\times \R}|x-y|^p\pi(dx,dy)$ over ${\cal M}(\mu,\check \mu^\Gamma)$ is the distribution of $(X,{\rm Proj}_{\Gamma}^{del}(X,U))$. Last, according to the remark after Proposition~10 in~\cite{PaWi1}, when $\mu\le_{cvx}\eta$ with $\eta$ compactly supported in $[x_1,x_N]$, then $\check\mu^\Gamma\le_{cvx}\check\eta^\Gamma$. This can be seen using the affine interpolation on $\Gamma$
$$\check\varphi^\Gamma(x):=1_{(-\infty,x_1)\cup [x_N,+\infty)}(x)\varphi(x)+\sum_{i=1}^{N-1}1_{[x_i,x_{i+1})}(x)\left(\frac{x_{i+1}-x}{x_{i+1}-x_i}\varphi(x_i)+\frac{x-x_{i}}{x_{i+1}-x_i}\varphi(x_{i+1})\right)$$
of a convex function $\varphi:\R\to\R$. Indeed $\check\varphi^\Gamma$ is still convex and one has
$$\int_\R\varphi(x)\check\mu^\Gamma(dx)=\int_\R\check\varphi^\Gamma(x)\mu(dx)\le \int_\R\check\varphi^\Gamma(x)\eta(dx)=\int_\R\varphi(x)\check\eta^\Gamma(dx).$$
According to the introduction in \cite{AlJo}, this convex order preservation does not generalize to higher dimensions where the minimizers are not so easy to express.

Whatever the dimension $d\in\N^*$, for every level $N\ge d+1$,  there exists an $L^p$-optimal dual quantization grid $\Gamma^{del}_{p,N}$ and a splitting operator ${\rm Proj}_{\Gamma^{del}_{p,N}}^{del}$ (see~\cite{PaWi1}) such that $$d_{p,N}(\mu)=d_p(\mu,\Gamma^{del}_{p,N})=\|X-{\rm Proj}_{\Gamma^{del}_{p,N}}^{del}(X,U)\|_p$$
and ${\rm Proj}_{\Gamma^{del}_{p,N}}^{del}(X,U)$ takes each value in $\Gamma^{del}_{p,N}$ with positive probability. For more details on this dual projection, see~\cite{PaWi1, PaWi2} where this notion has been developed and analyzed. We will see in the examples that even in dimension one, the convex order is not preserved by optimal dual quantization.

\smallskip
Notice that by Proposition~\ref{propprimcvx} $(b)$, the inequality $\cW_p^p(\mu,\nu)\le\cM_p^p(\mu,\nu)$ valid for $\nu\in{\cal P}_{\ge \mu}(\R^d,N)$ and \eqref{dpN}, \begin{equation}
   e_{p,N}(\mu)=W_p(\mu,{\cal P}(\R^d,N))\le W_p(\mu,{\cal P}_{\ge \mu}(\R^d,N))\le \inf_{\nu\in{\cal P}_{\ge \mu}(\R^d,N)}M_p(\nu,\mu)=d_{p,N}(\mu).\label{comprimdualproj}
\end{equation}
We may wonder whether the last inequality is an equality. Combining the tightness of any sequence of probability measures in ${\mathcal P}_{\ge \mu}(\R^d,N)$ minimizing the $W_p$-distance to $\mu$ deduced from the inequality $$\forall \nu\in{\mathcal P}(\R^d),\;\int_{\R^d}|x|^p\nu(dx)\le 2^{p-1}\left(\int_{\R^d}|x|^p\mu(dx)+W_p^p(\mu,\nu)\right),$$the closedness of ${\mathcal P}_{\ge \mu}(\R^d,N)$ for the weak convergence topology and the lower semi-continuity of the Wasserstein distance for this topology (see for instance Remark 6.12 p97 \cite{Villani}), we obtain the existence of a $W_p$-projection $\tilde\mu$ of $\mu$ on ${\mathcal P}_{\ge \mu}(\R^d,N)$. The last inequality in \eqref{comprimdualproj} is an equality iff the set
$$\left\{\pi\in{\cal P}(\mu,\tilde\mu):\cW^p_p(\mu,\tilde \mu)=\int_{\R^d\times\R^d}|y-x|^p \pi(dx,dy)\right\}$$
of $\cW_p$-optimal couplings between $\mu$ and  some $W_p$-projection $\tilde \mu$ of $\mu$ on ${\mathcal P}_{\ge \mu}(\R^d,N)$ intersects ${\mathcal M}(\mu,\tilde \mu)$. Moreover, $\tilde\mu$ is then the distribution of an $L^p$-optimal dual $N$-quantization of $\mu$. But there is no reason why the intersection should be non empty. We also may wonder whether, by a somewhat naive symmetry with the situation described in Proposition~\ref{propprimcvx} $(b)$ for the Voronoi quantization, the distribution of an $L^p$-optimal dual $N$-quantization of $\mu$ coincides with a $\cW_p$-projection $\tilde \mu$ of $\mu$ on ${\mathcal P}_{\ge \mu}(\R^d,N)$. According to the examples below, this property holds when $\mu$ is the uniform distribution on the interval $[0,1]$ (nevertheless $W_p({\mathcal U}[0,1],{\cal P}_{\ge {\mathcal U}[0,1]}(\R^d,N))<d_{p,N}({\mathcal U}[0,1])$) but is not true in general.

Note that, in the quadratic case $p=2$, for $\nu\in{\cal P}_{\ge \mu}(\R^d,N) $, since $M_2^2(\mu,\nu)=\int_{\R^d}|y|^2\nu(dy)-\int_{\R^d}|x|^2\mu(dx)$,
\begin{equation}
   d_{2,N}(\mu)^2=\inf_{\nu\in{\mathcal P}_{\ge \mu}(\R^d,N)}\int_{\R^d}|y|^2\nu(dy)-\int_{\R^d}|x|^2\mu(dx).\label{eqd22N}
\end{equation}
\begin{Proposition}
  In the quadratic case $p=2$, any optimal dual quantization grid $\Gamma^{del}_{2,N}$ remains optimal for each probability measure $\eta$ such that $\mu\le_{cvx}\eta\le_{cvx}\check\mu_{N}$ where $\check\mu_{N}$ denotes the distribution of ${\rm Proj}_{\Gamma^{del}_{2,N}}^{del}(X,U)$. 
  \end{Proposition}
 
 \noindent {\bf Proof.}
 Let $\eta$ be such that $\mu\le_{cvx}\eta\le_{cvx}\check\mu_N$. Since $\inf_{\nu\in{\mathcal P}_{\ge \mu}(\R^d,N)}\int_{\R^d}|y|^2\nu(dy)$ is attained for $\nu=\check\mu_N$, so is $\inf_{\nu\in{\mathcal P}_{\ge \eta}(\R^d,N)}\int_{\R^d}|y|^2\nu(dy)$. Therefore, \eqref{eqd22N} written with $\eta$ replacing $\mu$ implies that
 $$
d_{2,N}(\eta)^2=\int_{\R^d}|y|^2\check\mu_N(dy)-\int_{\R^d}|x|^2\eta(dx)=M_2^2(\eta,\check\mu_N).
$$
With the definition of $d_2(\eta,\Gamma^{del}_{2,N})$, we deduce that $d_{2,N}(\eta)^2\ge d_2(\eta,\Gamma^{del}_{2,N})^2$. Therefore $\Gamma^{del}_{2,N}$ is an optimal dual quadratic quantization grid for $\eta$. \hfill$\Box$
\bigskip

\noindent {\bf Examples.} $(a)$ Let $\mu={\cal U}[0,1]$ where, for two real numbers $a<b$,  ${\cal U}[a,b]$ denotes the uniform distribution on $[a,b]$ with density $\frac{1_{[a,b]}(x)}{b-a}$ with respect to the Lebesgue measure. We consider its approximation by probability measures in ${\cal P}(\R,N)$. A generic element of ${\cal P}(\R,N)$ writes $$\nu_N=\sum_{k=1}^Np_k\delta_{x_k}$$ with $x_1\le x_2\le \hdots\le x_N$ and $(p_1,\hdots,p_N)\in[0,1]^N$ satisfying $\sum_{k=1}^N p_k=1$. We will consider the particular choices $\hat\mu_N=\frac{1}{N}\sum_{k=1}^N\delta_{\frac{2k-1}{2N}}$ and $\check\mu_N=\frac{1}{2(N-1)}\delta_{0}+\frac{1}{N-1}\sum_{k=2}^{N-1}\delta_{\frac{k-1}{N-1}}+\frac{1}{2(N-1)}\delta_1$ of the respective distributions of the optimal primal and dual quantizations of $\mu$ on $N$ points. For $p\ge 1$, let us recover that $\hat\mu_N$ is the $W_p$-projection of ${\cal U}[0,1]$ on ${\cal P}(\R,N)$ (consequence of Proposition~\ref{propprimcvx} $(b)$) and check that $\check\mu_N$ is the $W_p$-projection of ${\cal U}[0,1]$ on ${\cal P}_{\ge {\cal U}[0,1]}(\R,N)$. The image of $\mu$ by $(0,1)\ni u\mapsto F_{\nu_N}^{-1}(u)-u$ is equal to $$\eta_N:=p_1{\cal U}[x_1-p_1,x_1]+p_2{\cal U}[x_2-(p_1+p_2),x_2-p_1]+\hdots+ p_N{\cal U}[x_N-1,x_N-(p_1+\hdots+p_{N-1})].$$ 
When $\nu_N=\hat\mu_N$ (resp. $\nu_N=\check\mu_N$) then $\eta_N=\hat\eta_N:={\cal U}[-\frac{1}{2N},\frac{1}{2N}]$ (resp. $\eta_N=\check\eta_N:={\cal U}[-\frac{1}{2(N-1)},\frac{1}{2(N-1)}]$). Since $\hat\eta_N=N1_{[-\frac{1}{2N},\frac{1}{2N}]}(x)dx$ and $\eta_N$ has a density with respect to the Lebesgue measure with values in $\{0,1,2,\hdots,N\}$, $\left(\eta_N-\hat\eta_N\right)^+$ is supported on the complement of $[\frac{-1}{2N},\frac{1}{2N}]$ where $\left(\hat\eta_N-\eta_N\right)^+$ is supported. Since both measures share the same mass, we deduce that for $p\ge 1$,
$$\int_{\R}|x|^p\left(\eta_N-\hat\eta_N\right)^+(dx)\ge \int_\R|x|^p\left(\hat\eta_N-\eta_N\right)^+(dx)\mbox{ i.e.}\int_\R|x|^p\eta_N(dx)\ge \int_\R|x|^p\hat\eta_N(dx).$$
Using that in dimension $d=1$, the  comonotonous coupling is $W_p$-optimal, we conclude that
$$\forall \nu_N\in{\cal P}(\R,N),\;W_p^p({\cal U}[0,1],\nu_N)=\int_{\R}|x|^p\eta_N(dx)\ge \int_\R|x|^p\hat\eta_N(dx)=W_p^p({\cal U}[0,1],\hat \nu_N),$$
with strict inequality unless $\eta_N=\hat\eta_N$. When $\eta_N$ is supported on $[-\frac{1}{2N},\frac{1}{2N}]$ (which is clearly equivalent to $\eta_N=\hat\eta_N$), then for each $k\in\{1,\hdots,N\}$, $\frac{-1}{2N}\le x_k-(p_1+\hdots+p_k)$ and $x_k-(p_1+\hdots+p_{k-1})\le\frac{1}{2N}$ so that $p_k\le\frac 1N$. With the normalisation, we deduce that $p_\ell=\frac 1N$ for each $\ell\in\{1,\hdots,N\}$, which plugged in the two last inequalities implies that $x_k=\frac{2k-1}{2N}$ for each $k\in\{1,\hdots,N\}$ so that $\nu_N=\hat\mu_N$. Therefore, for all $p\ge 1$, $\hat\mu_N$ is the $W_p$-projection of ${\cal U}[0,1]$ on ${\cal P}(\R,N)$ and
$$W_p({\cal U}[0,1],{\cal P}(\R,N))=W_p({\cal U}[0,1],\hat\mu_N)=\left(\int_{\R}|x|^p\hat\eta_N(dx)\right)^{1/p}=\frac{1}{2(p+1)^{1/p}N}.$$

If $\nu_N\ge_{cvx}{\cal U}[0,1]$, then $x_1\le 0$ and $x_N\ge 1$ so that $[x_1-p_1,x_1]\cap [x_N-1,x_N-(p_1+\hdots+p_{N-1})]\subset\{0\}$ and $\eta_N$ has a density with respect to the Lebesgue measure with values in $\{0,1,\hdots,N-1\}$. By repeating the previous argument with $\hat\eta_N$ replaced by $\check\eta_N$, we obtain that
\begin{align*}
   \forall \nu_N\in{\cal P}_{\ge {\cal U}[0,1]}(\R,N),\;W_p^p({\cal U}[0,1],\nu_N)=\int_{\R}|x|^p\eta_N(dx)\ge \int_{\R}|x|^p\check\eta_N(dx)
  =W_p^p({\cal U}[0,1],\check\mu_N),\end{align*}
with strict inequality unless $\eta_N=\check\eta_N$. When $\nu_N\ge_{cvx}{\cal U}[0,1]$ and $\eta_N$ is supported on $[-\frac{1}{2(N-1)},\frac{1}{2(N-1)}]$ (which is clearly equivalent to $\eta_N=\check\eta_N$), then for each $k\in\{1,\hdots,N\}$, $\frac{-1}{2(N-1)}\le x_k-(p_1+\hdots+p_k)$ and $x_k-(p_1+\hdots+p_{k-1})\le\frac{1}{2(N-1)}$ so that $p_k\le\frac {1}{N-1}$ and with the reinforced inequalities $x_1\le 0$ and $x_N\ge 1$ due to the convex order, $p_1\le \frac{1}{2(N-1)}$ and $p_N\le\frac{1}{2(N-1)}$. With the normalisation, we deduce that $p_1=p_N=\frac{1}{2(N-1)}$ and $p_\ell=\frac{1}{N-1}$ for each $\ell\in\{2,\hdots,N-1\}$, which plugged in the two (reinforced) inequalities implies that $x_1=0$, $x_k=\frac{k-1}{N-1}$ for each $k\in\{2,\hdots,N-1\}$ and $x_1=1$ so that $\nu_N=\check\mu_N$. Therefore, for all $p\ge 1$, $\check\mu_N$ is the $W_p$-projection of ${\cal U}[0,1]$ on ${\cal P}_{\ge {\cal U}[0,1]}(\R,N)$ and
\begin{align*}
   W_p({\cal U}[0,1],{\cal P}_{\ge {\cal U}[0,1]}(\R,N))&=W_p({\cal U}[0,1],\check\mu_N)=\left(\int_{\R}|x|^p\check\eta_N(dx)\right)^{1/p}=\frac{1}{2(p+1)^{1/p}(N-1)}\\&<\left(\frac{2}{(p+1)(p+2)}\right)^{1/p}\frac{1}{N-1}=d_{2,N}({\mathcal U}[0,1]).
\end{align*}
Notice that when $\nu_N\ge_{cvx}{\cal U}[0,1]$, then the measures $(\eta_N-\check\eta_N)^+$ supported on $[-\frac{1}{2(N-1)},\frac{1}{2(N-1)}]^c$ and $(\check \eta_N-\eta_N)^+$ supported on $[-\frac{1}{2(N-1)},\frac{1}{2(N-1)}]$ share the same mass and barycenter so that for each convex function $\varphi:\R\to\R$, $\int_\R\varphi(x)(\eta_N-\check\eta_N)^+(dx)\ge\int_\R\varphi(x)(\check \eta_N-\eta_N)^+(dx)$ and $\eta_N\ge_{cvx}\check\eta_N$.

The probability distribution  ${\cal U}[0,1]$ is smaller in the convex order than its dual quantization $\mu_6:=\frac 15\delta_0+\frac{7}{30}\delta_{\frac{2}{5}}+\frac{1}{15}\delta_{\frac{7}{15}}+\frac{1}{15}\delta_{\frac{8}{15}}+\frac{7}{30}\delta_{\frac{3}{5}}+\frac 15\delta_1$ on the grid $\left\{0,\frac{2}{5},\frac{7}{15},\frac{8}{15},\frac{3}{5},1\right\}$. Note that $\mu_6$ is not comparable with $\check\mu_6=\frac{1}{10}\delta_0+\frac{1}{5}\left(\delta_{\frac{1}{5}}+\delta_{\frac{2}{5}}+\delta_{\frac{3}{5}}+\delta_{\frac{4}{5}}\right)+\frac{1}{10}\delta_1$ for the convex order. Indeed, we have
\begin{align*}
  &\mu_6=\frac{2}{5}\left(\frac 12\delta_0+\frac 12\delta_{\frac 25}\right)+\frac 15\left(\frac 16\delta_{\frac 25}+\frac{1}{3}\delta_{\frac{7}{15}}+\frac{1}{3}\delta_{\frac{8}{15}}+\frac 16\delta_{\frac 35}\right)+\frac{2}{5}\left(\frac 12\delta_{\frac 35}+\frac 12\delta_{1}\right)\\
  &\check \mu_6=\frac{2}{5}\left(\frac 14\delta_0+\frac 12\delta_{\frac 15}+\frac 14\delta_{\frac 25}\right)+\frac 15\left(\frac 12\delta_{\frac 25}+\frac 12\delta_{\frac 35}\right)+\frac{2}{5}\left(\frac 14\delta_{\frac 35}+\frac 12\delta_{\frac 45}+\frac 14\delta_{1}\right)\\
  &\mbox{ with }\frac 12\delta_0+\frac 12\delta_{\frac 25}\ge_{cvx}\frac 14\delta_0+\frac 12\delta_{\frac 15}+\frac 14\delta_{\frac 25},\;\frac 16\delta_{\frac 25}+\frac{1}{3}\delta_{\frac{7}{15}}+\frac{1}{3}\delta_{\frac{8}{15}}+\frac 16\delta_{\frac 35}\le_{cvx}\frac 12\delta_{\frac 25}+\frac 12\delta_{\frac 35}\\&\mbox{ and }\frac 12\delta_{\frac 35}+\frac 12\delta_{1}\ge_{cvx}\frac 14\delta_{\frac 35}+\frac 12\delta_{\frac 45}+\frac 14\delta_{1}
\end{align*}
so that $\int_{\R}\varphi(x)\mu_6(dx)>\int_{\R}\varphi(x)\check\mu_6(dx)$ (resp. $\int_{\R}\varphi(x)\mu_6(dx)<\int_{\R}\varphi(x)\check\mu_6(dx)$) when the convex function $\varphi:\R\to\R$ is strictly convex on $[0,\frac{2}{5}]$ and affine on $[\frac{2}{5},1]$ (resp. affine on $[0,\frac 25]$, strictly convex on $[\frac 25,\frac 35]$ and affine on $[\frac 35,1]$). Since $\mu_6$ is clearly equal to its $L^p$-optimal dual $6$-quantization, this shows that the convex order is not preserved by optimal dual quantization.

\smallskip
\noindent $(b)$ Let now $\mu(dx)=2x\mbox{\bf 1}_{[0,1]}(x)dx$. We look for $\nu\in{\mathcal P}_{\ge \mu}(\R,3)$ minimizing either $\int_{\R}y^2\nu(dy)$ to compute the law of the optimal quadratic dual quantization of $\mu$ on $N=3$ points or $\cW_2^2(\mu,\nu)$ to compute $\tilde\mu$. Since $d=1$, $\cW_2^2(\mu,\nu)$ is equal to the integral $\int_0^1(F_\mu^{-1}(u)-F_\nu^{-1}(u))^2du$ of the squared difference between the quantile functions of $\mu$ and $\nu$. For the first criterion, we are going to check that it is equivalent to minimize over the following parametric subset of ${\mathcal P}_{\ge \mu}(\R,3)$ $$\left\{\nu_u(dy)=\frac{u}{3}\delta_0(dy)+\frac{1+\sqrt{u}}{3}\delta_{\sqrt{u}}(dy)+\frac{2-\sqrt{u}-u}{3}\delta_1(dy):u\in(0,1)\right\}.$$
One has
$\int y^2\nu_u(dy)=\frac{2+u^{3/2}-\sqrt{u}}{3}$ and the infimum is attained for $u=1/3$ so that $\Gamma^{del}_{2,3}=\{0,\frac{1}{\sqrt{3}},1\}$, $d_{2,3}(\mu)^2=d_{2}(\mu,\Gamma^{del}_{2,3})^2=\int y^2\nu_{1/3}(dy)-\int x^2\mu(dx)=\frac{1}{6}-\frac{2}{3^{5/2}}$. On the other hand, 
\begin{align*}
  \cW_2^2(\mu,\nu_u)&=\int_0^{u/3}(0-\sqrt{v})^2dv+\int_{u/3}^{(1+\sqrt{u}+u)/3}(\sqrt{u}-\sqrt{v})^2dv+\int_{(1+\sqrt{u}+u)/3}^1(1-\sqrt{v})^2dv\\
  &=-\frac{1}{6}+\frac{u^{3/2}-\sqrt{u}}{3}+4\frac{(1-\sqrt{u})(1+\sqrt{u}+u)^{3/2}+u^2}{3^{5/2}}.
\end{align*}
One easily checks that $\frac{d}{du}\cW_2^2(\mu,\nu_u)|_{u=1/3}>0$ and that $W_2^2(\mu,\nu_u)$ is minimal for $u\simeq 0.326$ so that $\tilde \mu\neq \nu_{1/3}$. Moreover, since ${\cal P}_{\ge \mu}(\Gamma^{del}_{2,3})=\{v\nu_{1/3}+(1-v)\nu_0:v\in[0,1]\}$ contains $\nu_{1/3}$, $$\cW_2^2(\mu,{\cal P}_{\ge \mu}(\Gamma^{del}_{2,3}))\le\cW^2_2(\mu,\nu_{1/3})\simeq 0.0199758<0.0383666\simeq d_{2}(\mu,\Gamma^{del}_{2,3})^2.$$

Let us finally check that $\inf_{\nu\in{\mathcal P}_{\ge \mu}(\R,3)}\int_\R y^2\nu(dy)\ge \inf_{u\in(0,1)}\int_\R y^2\nu_u(dy)$ and that  $\nu_u\in{\mathcal P}_{\ge \mu}(\R,3)$ for each $u\in(0,1)$. First note that for each $u\in(0,1)$, the mean $2/3$ of $\nu_u$ is equal to the one of $\mu$.
According to the characterization of the convex order in terms of potential functions, we have $\nu\ge_{cvx}\mu$ iff $\int_\R x\nu(dx)=\frac{2}{3}$ and \begin{equation}
   \forall x\in\R,\;\varphi_\nu(x):=\int_{-\infty}^x\nu((-\infty,y])dy\ge\varphi_\mu(x)  =1_{[0,1]}(x)\frac{x^3}{3}+1_{\{x>1\}}\left(x-\frac{2}{3}\right).\label{compote}
 \end{equation}
If $\nu$ weights at most three points, then the convex function $\varphi_\nu$ is piecewise affine with at most three changes of slope, the left-most slope being equal to $0$ and the right-most equal to $1$.
Therefore if $\nu\in{\mathcal P}_{\ge \mu}(\R,3)$ then $\nu(\{a\})>0$ for some $a\le 0$ and $\nu(\{b\})>0$ for some $b\ge 1$. If, moreover, $\nu(\{\sqrt{u}\})>0$ for some $u\in(0,1)$, then since the slope of $\varphi_\nu$ is constant on $(0,\sqrt{u})$ and on $(\sqrt{u},1)$, for all $x\in\R$,  $$\varphi_\nu(x)\ge 1_{(0,\sqrt{u}]}(x) \frac{u^{3/2}}{3}\times\frac{x}{\sqrt{u}}+1_{(\sqrt{u},1]}(x)\left(\frac{u^{3/2}}{3}\times \frac{1-x}{1-\sqrt{u}}+\frac{1}{3}\times \frac{x-\sqrt{u}}{1-\sqrt{u}}\right)+1_{\{x>1\}}\left(x-\frac{2}{3}\right)=\varphi_{\nu_u}(x),$$
so that, by convexity of the square function, $\int_\R y^2 \nu(dy)\ge \int_\R y^2\nu_u(dy)$. If, on the other hand, $\nu((0,1))=0$, then $\varphi_\nu$ has a constant slope on $(0,1)$ and we even have $$\int_\R y^2 \nu(dy)\ge \sup_{u\in(0,1)}\int_\R y^2\nu_u(dy).$$ Therefore $\inf_{\nu\in{\mathcal P}_{\ge \mu}(\R,3)}\int_\R y^2\nu(dy)\ge \inf_{u\in(0,1)}\int_\R y^2\nu_u(dy)$. To conclude that this inequality is an equality, it is enough to check that $\nu_u\in {\mathcal P}_{\ge \mu}(\R,3)$ for each $u\in(0,1)$. This follows from the inequality $\varphi_{\nu_u}(x)\ge \varphi_\mu(x)$ valid for all $x\in\R$ and all $u\in[0,1]$ since the graph of the convex function $\varphi_\mu$ is under its chords.

We finally recall the main result on convergence rate of dual quantization for bounded random vectors established in~\cite{PaWi3}. 
\begin{Theorem}[Pierce Lemma for dual quantization]\label{thm:zadoretpierce2} 
% $(a)$ {\sc Zador's Theorem  for dual quantization}: Let $X\!\in L_{\R^d}^{\infty}(\Omega,{\cal A}, \P)$ be a bounded random vector with distribution $\P_{_X}= \varphi.\lambda_d\stackrel{\perp}{+}\nu_{_X}$ where $\lambda_d$ denotes the Lebesgue measure and $\nu_{_X}$  denotes its  singular component. Then, for every $p\!\in (0, +\infty)$, 
% \[
% \lim_{N\to+\infty} N^{\frac 1d} d_{p,N}(X) = \widetilde J^{del}_{d,p}\left(\int_{\R^d} \varphi^{\frac{d}{d+p}}d\lambda_d\right)^{\frac 1d +\frac 1p}
% \]
% where $\displaystyle \widetilde J^{del}_{d,p}=\inf_{N\ge 1} N^{\frac 1d} d_{p,N}\big(\mathcal{U}([0,1]^d)\big)\ge \widetilde J^{vor}_{d,p}$.
% %($J_{p,d}$ is the similar  constant from the original Zador's Theorem). \textcolor{red}{ATTENTION! Les notations $J_{p,d}$ c'est pour la puissance p ieme, je crois\dots)}
% %
% %\smallskip
% When $d=1$, $\widetilde J^{del}_{1,p}=  \Big(\frac{2}{(p+1)(p+2)}\Big)^{1/p}$. Hence, $\frac{\widetilde J^{del}_{1,p}}{\widetilde J^{vor}_{1,p}}= \left(\frac{2^{p+1}}{p+2}\right)^{1/p} \uparrow 2$ as $p\uparrow +\infty$.

% \medskip
% \noindent $(b)$ {\sc Non-asymptotic bound (Pierce lemma)}: 
Let $p\ge 1$ and $\eta >0$. For every dimension $d\ge 1$, there exists a real constant $\widetilde C^{del}_{d,\eta, p}>0$ such that, for every  random vector $X:(\Omega,{\cal A}, \P) \to \R^d$, $L^{\infty}(\P)$-bounded, 
\begin{equation}\label{eq:Zadornonasymp}
d_{p,N}(X)\le \widetilde C^{del}_{d,\eta,p} N^{-\frac 1d} \sigma_{p+\eta}(X)
\end{equation}
where, for every $r>0$, $\sigma_r(X)= \inf_{a\in \R^d}\|X-a\|_r <+\infty$.
\end{Theorem}

\noindent {\bf Remark.} Note that this claim and the one in Theorem~\ref{thm:zadoretpierce1} remain true if the support of $\P_{_X}$ does not span $\R^d$ as an affine space, but $A_{\mu}$ with dimension $d'<d$. However, if  such is the case, then \eqref{majopierceprim} and~\eqref{eq:Zadornonasymp} hold with factor  $N^{-1/d'}$ so that $N^{-1/d}$ is suboptimal.
\section{Quantized approximations of martingale couplings}\label{secstab}
Let $\nu\in{\mathcal P}(\R^d)$ be compactly supported and $\mu\in{\cal P}_2(\R^d)$. For $K,N\ge 1$, let $\check{\nu}^K$ be an $L^p$-optimal dual $K$-quantization of $\nu$ with grid $\Gamma^{del}_{p,K}$ and $\hat{\mu}^N$ be a quadratic optimal primal $N$-quantization of $\mu$ with grid $\Gamma_{2,N}$. The two quantization errors correspond to martingale quasi-metrics (in comparison to Wasserstein metrics \eqref{defwas}, only martingale couplings are considered in the minimization defining \eqref{defwasmart}) between $\nu$ (resp. $\hat\mu_N$) and $\check{\nu}^K$ 
(resp. $\mu$) :

\begin{equation}
   d_{p,K}^p(\nu)=M_p^p(\nu,\check\nu^K)\mbox{ and }e_{2,N}^2(\mu)=M_2^2(\hat \mu^N,\mu).\label{dismart}
 \end{equation}
In contrast with  dual quantization where the martingale quasi-metric appears from the very beginning of the construction, the optimization in primal quantization relies on Wasserstein metrics. But in the quadratic $p=2$ case, the stationarity property \eqref{eq:VoroStatio} satisfied at optimality implies that $W_2(\hat \mu_N,\mu)=M_2(\hat \mu_N,\mu)$. 

Let $U\sim{\mathcal U}[0,1]$ and $q_y(d\check y)$ denote the law of ${\rm Proj}^{del}_{\Gamma^{del}_{p,K}}(y,U)$ when $y\in{\rm Conv}(\Gamma^{del}_{p,K})$ and $\delta_y(d\check y)$ otherwise. When $\mu\le_{cvx}\nu$, we are now going to exploit \eqref{dismart} to approximate any $\pi\in{\mathcal M}(\mu,\nu)$ by 
\begin{equation}
  \bar\pi^{N,K}(d\hat x,d\check y)=\int_{(x,y)\in\R^d\times\R^d}\delta_{{\rm Proj}_{\Gamma_{2,N}}(x)}(d\hat x)\pi(dx,dy)q_y(d\check y).\label{approcmartcoup} 
\end{equation}
For $(X,Y)\sim\pi$ independent from the random variable $U$ uniformly distributed on $[0,1]$, the random vector $({\rm Proj}_{\Gamma_{2,N}}(X),X,Y,{\rm Proj}_{\Gamma^{del}_{p,K}}^{del}(Y,U))$ is distributed according to $\delta_{{\rm Proj}_{\Gamma_{2,N}}(x)}(d\hat x)\pi(dx,dy)q_y(d\check y)$ and therefore $({\rm Proj}_{\Gamma_{2,N}}(X),{\rm Proj}_{\Gamma^{del}_{p,K}}^{del}(Y,U))$ is distributed according to $\bar\pi^{N,K}$. Therefore, the first marginal of $\bar\pi^{N,K}$ is the distribution $\hat\mu^N$ of ${\rm Proj}_{\Gamma_{2,N}}(X)$ and its second marginal is the law $\check\nu^K$ of ${\rm Proj}_{\Gamma^{del}_{p,K}}^{del}(Y,U)$. Moreover, using that $U$ is independent from $(X,Y)$ and \eqref{eq:Statio1} for the second equality, then that $\pi$ is a martingale coupling for the fourth equality and the stationarity property \eqref{eq:VoroStatio} for the last one, we obtain
\begin{align*}
   \E[{\rm Proj}_{\Gamma^{del}_{p,K}}^{del}(Y,U))|{\rm Proj}_{\Gamma_{2,N}}(X)]&=\E[\E[{\rm Proj}_{\Gamma^{del}_{p,K}}^{del}(Y,U))|(X,Y)]|{\rm Proj}_{\Gamma_{2,N}}(X)]=\E[Y|{\rm Proj}_{\Gamma_{2,N}}(X)]\\&=\E[\E[Y|X]|{\rm Proj}_{\Gamma_{2,N}}(X)]=\E[X|{\rm Proj}_{\Gamma_{2,N}}(X)]={\rm Proj}_{\Gamma_{2,N}}(X).
\end{align*}
so that $\bar\pi^{N,K}\in{\cal M}(\hat\mu^N,\check\nu^K)$.

\begin{Theorem}\label{thmapproxcoupl}Let $p\ge 1$, $\mu,\nu\in{\mathcal P}(\R^d)$ be such that $\mu\lecx\nu$ with $\nu$ compactly supported and for $N,K\ge 1$,  let $\hat{\mu}^N$ be a quadratic optimal primal $N$-quantization of $\mu$ and $\check{\nu}^K$ be an $L^p$-optimal dual $K$-quantization of $\nu$.
  Then, for each $\pi\in{\mathcal M}(\mu,\nu)$, the couplings $\bar\pi^{N,K}\in{\mathcal M}(\hat{\mu}^N,\check{\nu}^K)$ defined in \eqref{approcmartcoup} are such that
\begin{equation*}
   \cW_p^p(\bar\pi^{N,K},\pi)\le \begin{cases}
     e^p_{2,N}(\mu)+d^p_{p,K}(\nu)\mbox{ if $p\le 2$}\\
     CN^{-p/d}+2^{\frac{(p-2)}{2}}d^p_{p,K}(\nu)\mbox{ for $C<\infty$ not depending on $N,K$ if $2<p<2+d$}
   \end{cases}.
 \end{equation*}
 Moreover, when $p\ge 2$, $$\cW_2^2(\bar\pi^{N,K},\pi)\le e^2_{2,N}(\mu)+d_{p,K}^2(\nu).$$
 Last, for any $p\ge 1$, $$\limsup_{N\to\infty}\cAW_p(\bar\pi^{N,K},\pi)\le d_{p,K}(\nu).$$
\end{Theorem}According to Theorems~\ref{thm:zadoretpierce1} and~\ref{thm:zadoretpierce2}, $\sup_{N\ge 1}N^{1/d}e_{2,N}(\mu)<\infty$ and $\sup_{K\ge 1}K^{1/d}d_{p,K}(\nu)<\infty$.

\medskip
\noindent {\bf Proof.}
We have
\begin{align*}
   \cW_p^p(\bar\pi^{N,K},\pi)&\le \int_{\R^d\times\R^d\times\R^d\times\R^d}\left(|\hat x-x|^2+|y-\check y|^2\right)^{p/2}\delta_{{\rm Proj}_{\Gamma_{2,N}}(x)}(d\hat x)\pi(dx,dy)q_y(d\check y)\\&\le 2^{\frac{(p-2)^+}{2}}\int_{\R^d\times\R^d\times\R^d\times\R^d}\left(|\hat x-x|^p+|y-\check y|^p\right)\delta_{{\rm Proj}_{\Gamma_{2,N}}(x)}(d\hat x)\pi(dx,dy)q_y(d\check y)\\&=2^{\frac{(p-2)^+}{2}}\int_{\R^d}|{\rm Proj}_{\Gamma_{2,N}}(x)-x|^p\mu(dx)+2^{\frac{(p-2)^+}{2}}\E\left[\left|Y-{\rm Proj}^{del}_{\Gamma^{del}_{p,K}}(Y,U)\right|^p\right]\\&=2^{\frac{(p-2)^+}{2}}\int_{\R^d}|{\rm Proj}_{\Gamma_{2,N}}(x)-x|^p\mu(dx)+2^{\frac{(p-2)^+}{2}}d_{p,K}^p(\nu).
\end{align*}
When $p\le 2$, $\int_{\R^d}|{\rm Proj}_{\Gamma_{2,N}}(x)-x|^p\mu(dx)\le \left(\int_{\R^d}|{\rm Proj}_{\Gamma_{2,N}}(x)-x|^2\mu(dx)\right)^{p/2}=e^p_{2,N}(\mu)$. On the other hand, when $2<p<2+d$, since $\mu$ is compactly supported, we may apply the $L^2-L^p$- distortion mismatch Theorem 4.3 \cite{PaSa3} to obtain that $\sup_{N\ge 1}N^{p/d}\int_{\R^d}|{\rm Proj}_{\Gamma_{2,N}}(x)-x|^p\mu(dx)<\infty$, which completes the proof of the first inequality.
In a similar way, when $p\ge 2$
\begin{align*}
 \cW_2^2(\bar\pi^{N,K},\pi)&\le \int_{\R^d\times\R^d\times\R^d\times\R^d}\left(|\hat x-x|^2+|y-\check y|^2\right)\delta_{{\rm Proj}_{\Gamma_{2,N}}(x)}(d\hat x)\mu(dx)\pi_x(dy)q_y(d\check y)\\&=e_{2,N}^2(\mu)+\E\left[\left|Y-{\rm Proj}^{del}_{\Gamma^{del}_{p,K}}(Y,U)\right|^2\right]\le e_{2,N}^2(\mu)+\E\left[\left|Y-{\rm Proj}^{del}_{\Gamma^{del}_{p,K}}(Y,U)\right|^p\right]^{2/p}\\&=e_{2,N}^2(\mu)+d_{p,N}^2(\nu).
\end{align*}
Let now $$\check\pi^K(dx,d\check y)=\int_{y\in\R^d}\pi(dx,dy)q_y(d\check y)=\int_{y\in\R^d}\mu(dx)\pi_x(dy)q_y(d\check y).$$
We have $\check\pi^K\in{\cal M}(\mu,\check\nu^K)$.
Using the identity coupling $\mu(dx)\delta_x(d\tilde x)$ between $\mu$ and $\mu$ in the definition of the adapted Wasserstein distance then the coupling $\pi_x(dy)q_y(d\check y)$ between $\check\pi^K_x(d\check y)$ and $\pi_x(dy)$ in the definition of the usual Wasserstein distance, we obtain that
 \begin{align*}
   \cAW^p_p(\check\pi^K,\pi)&\le\int_{\R^d}\cW^p_p(\check\pi^K_x,\pi_x)  \mu(dx)\le \int_{\R^d}\int_{\R^d}\int_0^1|{\rm Proj}^{del}_{\Gamma^{del}_{p,K}}(y,u)-y|^pdu\pi_x(dy) \mu(dx)\\
                          &\le \int_{\R^d}\int_0^1|y-{\rm Proj}^{del}_{\Gamma^{del}_{p,K}}(y,u)|^pdu\nu(dy)=\E|Y-{\rm Proj}^{del}_{\Gamma^{del}_{p,K}}(Y,U)|^p=d^p_{p,K}(\nu).
 \end{align*}
With the triangle inequality, we deduce that to check the last statement in the theorem, it is enough to prove that $\lim_{N\to\infty}\cAW_p(\bar\pi^{N,K},\check\pi^K)=0$.

To do so, we denote by $(x_i)_{1\le i\le |\Gamma_{2,N}|}$ the $|\Gamma_{2,N}|\le N$ points in $\Gamma_{2,N}$ and by $(y_j)_{1\le j\le |\Gamma^{del}_{p,K}|}$ the $|\Gamma^{del}_{p,K}|\le K$ points in $\Gamma^{del}_{p,K}$ (as soon as the support of $\mu$ (resp. $\nu$) is not restricted to less than $N$ (resp. $K$) points, $|\Gamma_{2,N}|=N$ (resp. $|\Gamma^{del}_{p,K}|=K)$). Let $C_i=\{x\in\R^d:|x-x_i|<\min_{1\le k\le|\Gamma_{2,N}|,k\neq i}|x-x_k|\}$, $i\in\{1,\hdots,|\Gamma_{2,N}|\}$, denote the open Voronoi cells induced by $\Gamma_{2,N}$. Using \eqref{munechargepasfront} for the equality, we have
\begin{align}
  \cAW^p_p(\bar\pi^{N,K},\check\pi^K)&\le \int_{\R^d\times\R^d}(|\hat x-x|^p+\cW_p^p(\bar\pi^{N,K}_{\hat x},\check\pi^K_x))\delta_{{\rm Proj}_{\Gamma_{2,N}}(x)}(d\hat x)\mu(dx)\notag\\
  &=\sum_{i=1}^{|\Gamma_{2,N}|}\int_{C_i}(|x_i-x|^p+\cW_p^p(\bar\pi^{N,K}_{x_i},\check\pi^K_x))\mu(dx).\label{majadapwas}
\end{align}
Since $\check{\nu}^K(\Gamma^{del}_{p,K})=1$, $\mu(dx)$ a.e. $\check{\pi}^K_x(d\check y)=\sum_{j=1}^{|\Gamma^{del}_{p,K}|}q_j(x)\delta_{y_j}(d\check y)$ for some measurable functions $q_j$ with values in $[0,1]$ and such that $\sum_{j=1}^{|\Gamma^{del}_{p,K}|}q_j=1$. Moreover, for $i\in\{1,\hdots,|\Gamma_{2,N}|\}$, $$\bar \pi^{N,K}_{x_i}(d\check y)=\frac{1}{\mu(C_i)}\sum_{j=1}^{|\Gamma^{del}_{p,K}|}\int_{C_i}q_j(\xi)\mu(d\xi)\delta_{y_j}(d\check y).$$
As a consequence, we have
\begin{align*}
 {\rm TV}(\bar \pi^{N,K}_{x_i},\check{\pi}^K_x)=\sum_{j=1}^{|\Gamma^{del}_{p,K}|}\left|q_j(x)-\frac{1}{\mu(C_i)}\int_{C_i}q_j(\xi)\mu(d\xi)\right|.
\end{align*}
With Theorem 6.15 p.103 \cite{Villani}, we deduce that
\begin{align*}
   \cW_p^p(\bar\pi^{N,K}_{x_i},\check\pi^K_x)=  2^{p-1}\max_{1\le j\le |\Gamma^{del}_{p,K}|}|y_j|^p{\rm TV}(\hat \pi^{N,K}_{x_i},\hat{\pi}^K_x)\le 2^{p-1}\max_{1\le j\le |\Gamma^{del}_{p,K}|}|y_j|^p\sum_{j=1}^{|\Gamma^{del}_{p,K}|}\left|q_j(x)-\frac{1}{\mu(C_i)}\int_{C_i}q_j(\xi)\mu(d\xi)\right|.
\end{align*} 
Plugging this inequality in \eqref{majadapwas}, we deduce that
\begin{align}
  \cAW^p_p(\bar\pi^{N,K},\check\pi^K)\le \sum_{i=1}^{|\Gamma_{2,N}|}\int_{C_i}|x_i-x|^p\mu(dx)+
  2^{p-1}\max_{1\le j\le |\Gamma^{del}_{p,K}|}|y_j|^p\sum_{j=1}^{|\Gamma^{del}_{p,K}|}\sum_{i=1}^{|\Gamma_{2,N}|}\int_{C_i} \left|q_j(x)-\frac{1}{\mu(C_i)}\int_{C_i}q_j(\xi)\mu(d\xi)\right|\mu(dx).\label{majoawp}
\end{align}
Using Jensen's inequality in the case $p\in[1,2]$, the $L^2-L^p$- distortion mismatch Theorem 4.3 \cite{PaSa3} when $2<p<2+d$ and H\"older's inequality when $p\ge 2+d$, we get 
\begin{equation*}
   \sum_{i=1}^{|\Gamma_{2,N}|}\int_{C_i}|x_i-x|^p\mu(dx)=\int_{\R^d}|{\rm Proj}_{\Gamma_{2,N}}(x)-x|^p\mu(dx)\le\begin{cases}
     e^{p}_{2,N}(\mu)\,\mbox{ if }p\in[1,2]\\CN^{-p/d}\;\mbox{ if }2<p<2+d\\
     2^{(p-3)\vee 0}\max_{1\le j\le |\Gamma^{del}_{p,K}|}|y_j|^{p-2}e^{2}_{2,N}(\mu)\\
     \hskip 4cm\mbox{ if }p\ge 2+d
   \end{cases},
 \end{equation*}
 where the constant $C$ does not depend on $N$. With Theorem~\ref{thm:zadoretpierce1}, we deduce that the first term in the right-hand side of \eqref{majoawp} goes to $0$ as $N\to\infty$.
  Since $\int_{\R^d}|q_j(x)|\mu(dx)\le 1$, by Theorem 3.14 p.69 \cite{Rudin}, there exists a sequence $(q^n_j)_{n\in\N}$ of continuous and compactly supported functions  on $\R^d$ such that $\lim_{n\to\infty}\int_{\R^d}|q_j^n(x)-q_j(x)|\mu(dx)=0$. Since $q_j$ takes its values in the interval $[0,1]$, we suppose that so do the functions $q_j^n$ up to replacing them by $0\vee q^n_j\wedge 1$ which affects neither the continuity and compact support property nor the convergence. We have
 \begin{align*}
   \bigg|\sum_{i=1}^{|\Gamma_{2,N}|}&\int_{C_i} \left|q_j(x)-\frac{1}{\mu(C_i)}\int_{C_i}q_j(\xi)\mu(d\xi)\right|\mu(dx)-\sum_{i=1}^{|\Gamma_{2,N}|}\int_{C_i} \left|q^n_j(x)-\frac{1}{\mu(C_i)}\int_{C_i}q^n_j(\xi)\mu(d\xi)\right|\mu(dx)\bigg|\notag\\
   &\le \sum_{i=1}^{|\Gamma_{2,N}|}\int_{C_i} \left|q_j(x)-q_j^n(x)-\frac{1}{\mu(C_i)}\int_{C_i}(q_j(\xi)-q_j^n(\xi))\mu(d\xi)\right|\mu(dx)\notag\\
  &\le \sum_{i=1}^{|\Gamma_{2,N}|}\int_{C_i} |q_j(x)-q_j^n(x)|\mu(dx)+\sum_{i=1}^{|\Gamma_{2,N}|}\left|\int_{C_i}(q_j(\xi)-q_j^n(\xi))\mu(d\xi)\right|\le 2\int_{\R^d}|q_j(x)-q_j^n(x)|\mu(dx),% 
 \end{align*}
 where the right-hand side goes to $0$ as $n\to\infty$.
We deduce that to prove that the second term in the right-hand side of \eqref{majoawp} goes to $0$ as $N\to\infty$, it is enough to check that so does $\sum_{i=1}^{|\Gamma_{2,N}|}\int_{C_i} \left|q_j(x)-\frac{1}{\mu(C_i)}\int_{C_i}q^n_j(\xi)\mu(d\xi)\right|\mu(dx)$ for any fixed $n\in\N$. For $X\sim\mu$ and $\hat X^N={\rm Proj}_{\Gamma_{2,N}}(X)$,

\begin{align*}
  \sum_{i=1}^{|\Gamma_{2,N}|}&\int_{C_i} \left|q^n_j(x)-\frac{1}{\mu(C_i)}\int_{C_i}q^n_j(\xi)\mu(d\xi)\right|\mu(dx)=\E\left|q^n_j(X)-\E[q^n_j(X)|\hat X^{N}]\right|\\&\le \E^{1/2}\left[\left(q^n_j(X)-\E[q^n_j(X)|\hat X^{N}]\right)^2\right]\le \E^{1/2}\left[\left(q^n_j(X)-q^n_j(\hat X^{N})\right)^2\right]
\end{align*}
where, for the last inequality, we used that the conditional expectation given $\hat X^{N}$ is the best quadratic approximation of a random variable by a measurable function of $\hat X^{N}$. Let $\varepsilon>0$. Since $q_j^n$ is continuous and compactly supported, this function is uniformly continuous. Since it takes its values in the interval $[0,1]$, we deduce that there exists $\eta>0$ such that for all $x,y\in\R^d$, $|q^n_j(x)-q^n_j(y)|\le \varepsilon 1_{\{|x-y|\le\eta\}}+1_{\{|x-y|>\eta\}}$. Therefore
\begin{align*}
   \E\left[\left(q^n_j(X)-q^n_j(\hat X^{N})\right)^2\right]\le \varepsilon^2+\P(|X-\hat X^{N}|\ge\eta)\le \varepsilon^2+\frac{\E[|X-\hat X^{N}|^2]}{\eta^2}=\varepsilon^2+\frac{e_2(\mu,N)^2}{\eta^2}.
\end{align*}
With Theorem~\ref{thm:zadoretpierce1}, we deduce that the left-hand side goes to $0$ as $N\to\infty$ and conclude that so does $\cAW_p(\bar\pi^{N,K},\check\pi^K)$.$\hfill\cqfd$

\section{Application to weak martingale optimal transport problems}

We endow ${\cal P}_{\le}\times {\cal P}_1(\R^d)=\left\{(\mu,\nu):\mu,\nu\in{\mathcal P}_1(\R^d)\mbox{ and }\mu\le_{cvx}\nu\right\}$ with the metric $W_1(\mu,\tilde \mu)+W_1(\nu,\tilde \nu)$ between $(\mu,\nu)$ and $(\tilde\mu,\tilde\nu)$ and $\R^d\times{\mathcal P}_1(\R^d)$ with the metric obtained as the sum of the Euclidean distance on $\R^d$ and the Wasserstein distance $W_1$ on ${\cal P}_1(\R^d)$. For a cost function $C:\R^d\times{\mathcal P}_1(\R^d)\to\R$ Borel measurable, the Weak Martingale Optimal Transport problem introduced in \cite{BackPam} consists in computing for $(\mu,\nu)\in{\cal P}_{\le}\times {\cal P}_1(\R^d)$
\begin{equation}
   V(\mu,\nu)=\inf_{\pi\in{\cal M}(\mu,\nu)}\int_{\R^d}C(x,\pi_x)\mu(dx)\label{mowt}
\end{equation}
and the minimal couplings $\pi\in{\cal M}(\mu,\nu)$. For the choice
\begin{equation}
   \tilde C(x,\eta)=\begin{cases}
     C(x,\eta)\mbox{ if }\int_{\R^d}y\eta(dy)=x\\
     +\infty\mbox{ otherwise}
   \end{cases},\label{ctildec}
\end{equation}it can be seen as a particular case of the Weak Optimal Transport problem \begin{equation}
   \tilde V(\mu,\nu)=\inf_{\pi\in{\cal P}(\mu,\nu)}\int_{\R^d}\tilde C(x,\pi_x)\mu(dx)\mbox{ for }(\mu,\nu)\in{\cal P}_1(\R^d)\times{\cal P}_1(\R^d)\label{owt}
\end{equation} introduced by Gozlan, Roberto, Samson and Tetali in \cite{GoRoSamTe} and studied by Backhoff-Veraguas, Beiglb\"ock and Pammer in \cite{BackBeiPam}. Indeed, when $(\mu,\nu)\in{\cal P}_{\le}\times {\cal P}_1(\R^d)$ then for each $\pi\in{\cal P}(\mu,\nu)\setminus{\cal M}(\mu,\nu)$, $\int_{\R^d}\tilde C(x,\pi_x)\mu(dx)=+\infty$ and for each $\pi\in{\cal M}(\mu,\nu)$, $\int_{\R^d}\tilde C(x,\pi_x)\mu(dx)=\int_{\R^d}C(x,\pi_x)\mu(dx)$, which implies that $\tilde V(\mu,\nu)=V(\mu,\nu)$.

The martingale optimal transport problem corresponds to the particular case of the WMOT problem when the cost function is linear in the measure component : $C(x,\eta)=\int_{\R^d} c(x,y)\eta(dy)$ for a Borel measurable function $c:\R^d\times\R^d\to\R$ with at most linear growth in its second variable.

The existence of minimal couplings in the WMOT problem \eqref{mowt} and the lower semi-continuity of the value function $V$ are deduced from Theorem 2.6 and Proposition 5.8 $(b)$ \cite{BJMP}. 
\begin{Proposition}\label{thmsci}
  Assume that $C:\R^d\times{\mathcal P}_1(\R^d)\to\R$ is lower semi-continuous, convex in the measure argument (for all $x\in\R^d$, ${\cal P}_1(\R^d)\ni\eta\mapsto C(x,\eta)$ is convex) and such that $$\sup_{(x,\eta)\in\R^d\times{\mathcal P}_1(\R^d)}\frac{|C(x,\eta)|}{1+|x|+\int_{\R^d}|y|\eta(dy)}<+\infty.$$ Then for each $(\mu,\nu)\in{\cal P}_{\le}\times {\cal P}_1(\R^d)$, there exists $\pi^\star\in{\cal M}(\mu,\nu)$, unique if $C$ is strictly convex in the measure argument, such that $V(\mu,\nu)=\int_{\R^d}C(x,\pi^\star_x)\mu(dx)$ and $(\mu,\nu)\mapsto V(\mu,\nu)$ is lower semi-continuous on ${\cal P}_{\le}\times {\cal P}_1(\R^d)$.
\end{Proposition}

Theorem 2.6 \cite{BJMP} also ensures convergence of the optimal couplings under convergence of the value function, a property which holds in dimension $d=1$.
\begin{Proposition}\label{propclsc}
  Assume that $C:\R^d\times{\mathcal P}_1(\R^d)\to\R$ is lower semi-continuous, convex in the measure argument and such that $$\sup_{(x,\eta)\in\R^d\times{\mathcal P}_1(\R^d)}\frac{|C(x,\eta)|}{1+|x|+\int_{\R^d}|y|\eta(dy)}<+\infty.$$
  Let for each $k\in\N$, $(\mu_k,\nu_k)\in{\cal P}_{\le}\times {\cal P}_1(\R^d)$ and $\pi^\star_k$ be an optimal coupling for $V(\mu_k,\nu_k)$, the existence of which is a consequence of Theorem~\ref{thmsci}.  If $((\mu_k,\nu_k))_{k\in\N}$ converges  to $(\mu,\nu)$ in ${\cal P}_{\le}\times {\cal P}_1(\R^d)$ as $k\to\infty$ and $V(\mu,\nu)=\lim_{k\to\infty}V(\mu_k,\nu_k)$, then all the accumulation points of $(\pi^\star_k)_{k\in\N}$ for the weak convergence topology are minimizers for $V(\mu,\nu)$. If $C$ is moreover strictly convex in the measure argument, then the sequence $(\pi^\star_k)_{k\in\N}$ converges in $AW_1$ to the unique optimal coupling $\pi^\star$ between $\mu$ and $\nu$.\\
  If $d=1$ and either $C$ is continuous or $C$ is continuous in its second argument and for each Borel subset $A$ of $\R^d$, $(\mu_k(A))_{k\in\N}$ converges to $\mu(A)$ as $k\to\infty$, then $V(\mu,\nu)=\lim_{k\to\infty}V(\mu_k,\nu_k)$.
 \end{Proposition}
When $(\mu_k,\nu_k)=(\hat{\mu}^{N_k},\check{\nu}^{K_k})$ with $\hat{\mu}^{N_k}$ a quadratic optimal primal $N_k$-quantization of $\mu$ and $\check{\nu}^{K_k}$ an $L^p$-optimal dual $K_k$-quantization of $\nu$, then Theorem~\ref{thmapproxcoupl} ensures that it is possible to approximate in $AW_1$ distance any optimal martingale coupling $\pi^\star$ between $\mu$ and $\nu$ by martingale couplings between $\hat{\mu}^{N_k}$ and $\check{\nu}^{K_k}$  and we deduce the upper-semicontinuity of the value function along this sequence, whatever the dimension $d$. 

\begin{Lemma}\label{lemccont}
  Let $p\ge 1$, $\mu,\nu\in{\mathcal P}(\R^d)$ be such that $\mu\lecx\nu$ with $\nu$ compactly supported and for $N,K\ge 1$,  $\hat{\mu}^N$ be a quadratic optimal primal $N$-quantization of $\mu$ and $\check{\nu}^{K}$ an $L^p$-optimal dual $K$-quantization of $\nu$.  
If $C:\R^d\times{\mathcal P}_1(\R^d)\to\R$ is continuous, then $V(\mu,\nu)$ is finite. If moreover, $x\mapsto \sup_{\eta\in{\cal P}_1(\R^d)}\frac{|C(x,\eta)|}{1+\int_{\R^d}|y|\eta(dy)}$ is locally bounded on $\R^d$, then $$\limsup_{N,K\to\infty}V(\hat{\mu}^N,\check\nu^K)\le V(\mu,\nu).$$
If $C:\R^d\times{\mathcal P}_1(\R^d)\to\R$ is Lipschitz continuous with constant ${\rm Lip}(C)$, then $$\limsup_{N\to\infty}V(\hat{\mu}^N,\check\nu^K)\le V(\mu,\nu)+{\rm Lip}(C)d_{p,K}(\nu).$$\end{Lemma}
With Propositions~\ref{thmsci} and~\ref{propclsc}, we easily deduce the following corollary.

\begin{Corollary}Let $p\ge 1$, $\mu,\nu\in{\mathcal P}(\R^d)$ be such that $\mu\lecx\nu$ with $\nu$ compactly supported and for $N,K\ge 1$,  $\hat{\mu}^N$ be a quadratic optimal primal $N$-quantization of $\mu$ and $\check{\nu}^K$ an $L^p$-optimal dual $K$-quantization of $\nu$. If $C:\R^d\times{\mathcal P}_1(\R^d)\to\R$ is continuous, convex in the measure argument and such that $$\sup_{(x,\eta)\in\R^d\times{\mathcal P}_1(\R^d)}\frac{|C(x,\eta)|}{1+|x|+\int_{\R^d}|y|\eta(dy)}<+\infty,$$then $\lim_{N,K\to\infty}V(\hat\mu^N,\check\nu^K)=V(\mu,\nu)$. Moreover, all accumulation points as $N,K\to\infty$ of sequences $(\pi^\star_{N,K})_{N,K}$ of minimizers for $V(\hat\mu^N,\check\nu^K)$ are minimizers for $V(\mu,\nu)$. If $C$ is also strictly convex in the measure argument, then $\lim_{N,K\to\infty}AW_1(\pi^\star_{N,K},\pi^\star)=0$, where $\pi^\star$ is the unique optimal coupling between $\mu$ and $\nu$.
\end{Corollary}

\noindent {\bf Proof of Lemma~\ref{lemccont}.} We first suppose that $C:\R^d\times{\mathcal P}_1(\R^d)\to\R$ is continuous. Let $\bar B$ be some closed ball centered at the origin with sufficiently large radius so that $\nu(\bar B)=1$ and ${\cal P}_{\bar B}(\R^d)=\{\eta\in{\cal P}(\R^d):\eta(\bar B)=1\}$. The set ${\cal P}_{\bar B}(\R^d)$ is a compact subset of ${\cal P}_1(\R^d)$ and therefore the continuous cost function $C$ is bounded on the compact subset $\bar B\times{\cal P}_{\bar B}(\R^d)$ of $\R^d\times{\cal P}_1(\R^d)$. Since $\hat\mu^N\le_{cvx}\mu\le_{cvx}\nu$, $\hat{\mu}^N(\bar B)=\mu(\bar B)=1$. Moreover for $\pi\in{\cal P}(\mu,\nu)$, $\int_{\R^d}\pi_x(\bar B)\mu(dx)=\nu(\bar B)=1$ so that $(x,\pi_x)\in\bar B\times{\cal P}_{\bar B}(\R^d)$ $\mu(dx)$ a.e.. Therefore $$\inf_{(x,\eta)\in\bar B\times{\cal P}_{\bar B}(\R^d)}C(x,\eta)\le V(\mu,\nu)\le \sup_{(x,\eta)\in\bar B\times{\cal P}_{\bar B}(\R^d)}C(x,\eta)$$
and $V(\mu,\nu)$ is finite. Let $\varepsilon>0$. By the continuity and the growth assumption satisfied by the cost function $C$ and the compactness of $\bar B\times{\cal P}_{\bar B}(\R^d)$ , 
\begin{align*}
  \exists \alpha>0,&\;\forall (x,\eta,\tilde x,\tilde \eta)\in\bar B\times{\cal P}_{\bar B}(\R^d)\times\R^d\times{\cal P}_1(\R^d)\mbox{ s.t. }|x-\tilde x|+{\cW}_1(\eta,\tilde \eta)\le \alpha,\;|C(x,\eta)-C(\tilde x,\tilde \eta)|\le \varepsilon\\
  &\mbox{ and }\forall (x,\eta,\tilde x,\tilde \eta)\in\bar{B}\times{\cal P}_{\bar B}(\R^d)\times \bar{B}\times{\cal P}_{1}(\R^d),\;|C(x,\eta)|+|C(\tilde x,\tilde \eta)|\le \frac{1}{\alpha}\left(1+\int_{\R^d}|y|\tilde\eta(dy)\right).
\end{align*} Let $\pi\in{\cal M}(\mu,\nu)$ be such that $\int_{\R^d}C(x,\pi_x)\mu(dx)\le V(\mu,\nu)+\varepsilon$. By Theorem~\ref{thmapproxcoupl}, there exists $\bar\pi^{N,K}\in{\cal M}(\hat{\mu}^N,\check\nu^K)$ such that $\limsup_{N\to\infty}AW_1(\bar\pi^{N,K},\pi)\le d_{p,K}(\nu)$ and $\lim_{N,K\to\infty}AW_1(\bar\pi^{N,K},\pi)=0$. Let $m^{N,K}\in{\cal P}(\mu,\hat\mu^N)$ be an optimal coupling for $AW_1(\pi,\bar\pi^{N,K})$. We have
\begin{align}
  &\left|\int_{\R^d}C(x,\pi_x)\mu(dx)-\int_{\R^d}C(\tilde x,\bar\pi^{N,K}_{\tilde x})\hat\mu^N(d\tilde x)\right|\le \int_{\R^d\times\R^d}\left|C(x,\pi_x)-C(\tilde x,\bar\pi^{N,K}_{\tilde x})\right|m^{N,K}(dx,d\tilde x)\notag\\
  &\le \varepsilon+\int_{\R^d\times\R^d}\left|C(x,\pi_x)-C(\tilde x,\bar\pi^{N,K}_{\tilde x})\right|1_{\{|x-\tilde x|+{\cW}_1(\pi_x,\bar\pi^{N,K}_{\tilde x})\ge \alpha\}}m^{N,K}(dx,d\tilde x)\notag\\
    &\le \varepsilon+\frac 1\alpha\int_{\R^d\times\R^d}\left(1+\int_{\R^d}|y|\bar\pi^{N,K}_{\tilde x}(dy)\right)1_{\{|x-\tilde x|+{\cW}_1(\pi_x,\bar\pi^{N,K}_{\tilde x})\ge \alpha\}}m^{N,K}(dx,d\tilde x).\label{majodiffcout}
\end{align}
By Markov inequality, $\int_{\R^d\times\R^d}1_{\{|x-\tilde x|+{\cW}_1(\pi_x,\bar\pi^{N,K}_{\tilde x})\ge \alpha\}}m^{N,K}(dx,d\tilde x)\le\frac{AW_1(\pi,\bar\pi^{N,K})}{\alpha}$. Moreover, since $$\int_{(x,\tilde x)\in\R^d\times\R^d}\bar\pi^{N,K}_{\tilde x}(dy)m^{N,K}(dx,d\tilde x)=\int_{\tilde x\in\R^d}\bar\pi^{N,K}_{\tilde x}(dy)\hat\mu^N(d\tilde x)=\int_{\tilde x\in\R^d}\bar\pi^{N,K}(d\tilde x,dy)=\check\nu^K(dy),$$
there is a Markov kernel $q_y(dx,d\tilde x)$ such that $\bar\pi^{N,K}_{\tilde x}(dy)m^{N,K}(dx,d\tilde x)=\check\nu^K(dy)q_y(dx,d\tilde x)$ and 
\begin{align*}
 \int_{\R^d\times\R^d}\int_{\R^d}|y|\bar\pi^{N,K}_{\tilde x}(dy)1_{\{|x-\tilde x|+{\cW}_1(\pi_x,\bar\pi^{N,K}_{\tilde x})\ge \alpha\}}m^{N,K}(dx,d\tilde x)=\int_{\R^d}|y|\beta(y)\check\nu^K(dy)
\end{align*}
where the function $\beta(y):=\int_{\R^d\times\R^d}  1_{\{|x-\tilde x|+{\cW}_1(\pi_x,\bar\pi^{N,K}_{\tilde x})\ge \alpha\}}q_y(dx,d\tilde x)$ is $[0,1]$-valued and such that $$\int_{\R^d}\beta(y)\check\nu^K(dy)=\int_{\R^d\times\R^d}1_{\{|x-\tilde x|+{\cW}_1(\pi_x,\bar\pi^{N,K}_{\tilde x})\ge \alpha\}}m^{N,K}(dx,d\tilde x)\le\frac{AW_1(\pi,\bar\pi^{N,K})}{\alpha}.$$
Since $\lim_{N,K\to\infty} AW_1(\pi,\bar\pi^{N,K})=0=\lim_{K\to\infty}W_1(\nu,\check\nu^K)$, with Lemma.. BJMP, we deduce that the second term in the right-hand side of \eqref{majodiffcout} tends to $0$ as $N,K\to\infty$. Hence
\begin{align*}
  \limsup_{N,K\to\infty} V(\hat\mu^N,\check\nu^K)&\le \limsup_{N,K\to\infty}\int_{\R^d}C(\tilde x,\bar\pi^{N,K}_{\tilde x})\hat\mu^N(d\tilde x)\\&\le \int_{\R^d}C(x,\pi_x)\mu(dx)+ \limsup_{N,K\to\infty}\left|\int_{\R^d}C(x,\pi_x)\mu(dx)-\int_{\R^d}C(\tilde x,\bar\pi^{N,K}_{\tilde x})\hat\mu^N(d\tilde x)\right|\\&\le V(\mu,\nu)+2\varepsilon. 
\end{align*}
Since $\varepsilon$ is arbitrary, we conclude that $\limsup_{N,K\to\infty} V(\hat\mu^N,\check\nu^K)\le V(\mu,\nu)$.

We now suppose that $C:\R^d\times{\mathcal P}_1(\R^d)\to\R$ is Lipschitz continuous with constant ${\rm Lip}(C)$. We then have
\begin{align*}
   &\left|\int_{\R^d}C(x,\pi_x)\mu(dx)-\int_{\R^d}C(\tilde x,\bar\pi^{N,K}_{\tilde x})\hat\mu^N(d\tilde x)\right|\le \int_{\R^d\times\R^d}\left|C(x,\pi_x)-C(\tilde x,\bar\pi^{N,K}_{\tilde x})\right|m^{N,K}(dx,d\tilde x)\\&\le {\rm Lip}(C)\int_{\R^d\times\R^d}\left(|x-\tilde x|+W_1(\pi_x,\bar\pi^{N,K}_{\tilde x})\right)m(dx,d\tilde x)={\rm Lip}(C)AW_1(\pi,\bar\pi^{N,K}).
\end{align*}
Therefore
\begin{align*}
   V(\hat\mu^N,\check\nu^K)\le \int_{\R^d}C(\tilde x,\bar\pi^{N,K}_{\tilde x})\hat\mu^N(d\tilde x)&\le \int_{\R^d}C(x,\pi_x)\mu(dx)+{\rm Lip}(C)AW_1(\pi,\bar\pi^{N,K})\\&\le V(\mu,\nu)+\varepsilon+{\rm Lip}(C)AW_1(\pi,\bar\pi^{N,K}).
\end{align*}
Since $\varepsilon$ is arbitrary, we deduce that $\limsup_{N\to\infty}V(\hat\mu^N,\check\nu^K)\le V(\mu,\nu)+{\rm Lip}(C)d_{p,K}(\nu)$.
$\hfill\cqfd$

\end{document}